\documentclass[11pt]{amsart}

\usepackage{amsmath}

\usepackage{amssymb}

\newcommand{\R}{{\mathbb R}}       
\newcommand{\Z}{{\mathbb Z}}       

\newcommand{\CC}{{\mathcal C}}

\newcommand{\dist}{{\rm dist}}
\newcommand{\ds}{\displaystyle }

\newcommand{\lra}{{\longrightarrow}}

\newcommand{\noi}{\noindent}

\newcommand{\rf}[1]{{(\ref{#1})}}

\newcommand{\supp}{{\rm supp}}

\newcommand{\vphi}{{\varphi}}
\newcommand{\ve}{{\varepsilon}}
\newcommand{\vv}{{\vspace{2mm}}}

\newcommand{\wt}[1]{{\widetilde{#1}}}
\newcommand{\wh}[1]{{\widehat{#1}}}

\newcommand{\bmo}{{B\!M\!O}}
\newcommand{\rbmo}{{R\!B\!M\!O}}

\newcommand{\QH}{{\wh{\wh{Q}^3 \hspace{-1mm}} }}

\newtheorem{theorem}{Theorem}[section]
\newtheorem{lemma}[theorem]{Lemma}
\newtheorem{coro}[theorem]{Corollary}

\theoremstyle{definition}
\newtheorem{definition}[theorem]{Definition}
\newtheorem{example}[theorem]{Example}

\theoremstyle{remark}
\newtheorem{rem}[theorem]{Remark}

\numberwithin{equation}{section}

\newcommand{\brem}{\begin{rem}}
\newcommand{\erem}{\end{rem}}

\newcommand{\bexam}{\begin{example}}
\newcommand{\eexam}{\end{example}}

\begin{document}

\title[Non doubling measures]
{Littlewood-Paley theory and the $T(1)$ theorem with non doubling measures}

\author[XAVIER TOLSA]{Xavier Tolsa}

\address{Department of Mathematics,
University of G\"oteborg / Chalmers, 412 96 G\"oteborg, Sweden}

\email{xavier@math.chalmers.se}

\thanks{Supported by a postdoctoral grant from the European Commission for the
TMR Network ``Harmonic Analysis''. Also
partially supported by grants DGICYT PB96-1183 and CIRIT
1998-SGR00052 (Spain).}


\date{May 21, 2000.}


\begin{abstract}
Let $\mu$ be a Radon measure on $\R^d$ which may be non doubling.
The only condition that $\mu$ must satisfy is
$\mu(B(x,r))\leq C\,r^n$, for all $x\in\R^d$, $r>0$ and for some fixed
$0<n\leq d$. In this paper, Littlewood-Paley theory for functions
in $L^p(\mu)$ is developed. One of the main difficulties to be solved is the
construction of ``reasonable'' approximations of the identity in order to
obtain a Calder\'on type reproducing formula.
Moreover, it is shown that the $T(1)$ theorem for $n$-dimensional
Calder\'on-Zygmund operators, without doubling assumptions, can be proved
using the Littlewood-Paley type decomposition that is obtained
for functions in $L^2(\mu)$, as in the classical case of homogeneous spaces.
\end{abstract}

\maketitle


\section{Introduction}

A basic hypothesis in the classical Calder\'on-Zygmund theory of harmonic
analysis is the doubling property of the underlying measure $\mu$ on $\R^d$
(or on more general spaces, such as the so called homogeneous spaces).
A measure $\mu$ is said to be doubling if there exists some constant
$C$ such that $\mu(B(x,2r))\leq C\,\mu(B(x,r))$ for all $x\in\supp(\mu)$,
$r>0$.
Recently it has been shown many results of the theory also hold
without assuming the doubling property. Some of these results,
such as the ones in \cite{NTV1}, \cite{NTV2},
\cite{NTV3}, \cite{Tolsa1}, \cite{Tolsa2}, deal with
Calder\'on-Zygmund operators.
Other questions are related to the
spaces $\bmo$ and $H^1$ \cite{MMNO}, \cite{Tolsa3},
\cite{Tolsa5}; or with vector valued inequalities and weights \cite{GM},
\cite{OP}, etc.

The aim of this paper is twofold. The first objective consists of
developing some Littlewood-Paley theory for functions in $L^p(\mu)$,
$1<p<\infty$, with $\mu$ being a Radon measure on $\R^d$ which may be non
doubling. The
only condition that $\mu$ must satisfy is the growth condition
\begin{equation}  \label{creix}
\mu(B(x,r))\leq C_0\,r^n \qquad\mbox{for all $x\in \supp(\mu)$, $r>0$,}
\end{equation}
where $n$ is some fixed number such that $0<n\leq d$.

The second objective of the paper is to apply these Littlewood-Paley
techniques to obtain a new proof of the $T(1)$ theorem for non doubling
measures
on $\R^d$ (see Theorem \ref{t1} below for the precise statement of the result).
The classical $T(1)$ theorem (with $\mu$ being the Lebesgue measure on $\R^d$)
was proved by David and Journ\'e \cite{DJ}. This result was extended recently
to the case of non doubling measures on $\R^d$ by Nazarov, Treil
and Volberg using dyadic martingales associated with random dyadic lattices.
Another proof in the setting of non doubling measures suitable for the Cauchy
integral operator was obtained at the same time independently by the author
\cite{Tolsa1}.
The proof of the $T(1)$ theorem that we will give in this paper will
follow an approach similar to the one of Coifman for proving this result
in the case of homogeneous spaces
(cf. \cite{DJS}), and to the one of David, Journ\'e and Semmes \cite{DJS} for
obtaining the $T(b)$ theorem for homogeneous spaces.

Let us remark that in the particular case of the Cauchy integral operator
other proofs of the $T(1)$ theorem have been given (see \cite{Verdera},
\cite{Tolsa3}) but,
as far as we know, for general Calder\'on-Zygmund operators the only proof
available for the moment was the one of Nazarov, Treil and Volberg based on
random dyadic lattices.

One of the main difficulties for developing Littlewood-Paley theory
with respect to some measure $\mu$ which does not satisfy any regularity
property, apart from the growth condition \rf{creix}, is the construction
of ``reasonable'' approximations of the identity.
Our geometric construction will
be based on some ideas originated from \cite{Tolsa5}, where an atomic
Hardy space useful for studying the $L^p(\mu)$ boundedness of
Calder\'on-Zygmund operators (with $\mu$ non doubling) was characterized in
terms of some grand maximal operator. A necessary step for the proof was the
construction of a suitable lattice of cubes and of smooth functions
$\vphi_{y,k}(x)$ associated to the corresponding cubes.
In the present paper we will use a
slight variant of this lattice. Moreover, the functions $\vphi_{y,k}(x)$
will play an essential role in our construction of the approximation
of the identity.

Once we have at our disposal this approximation of the identity, we will apply
some ideas of Coifman for obtaining a Calder\'on type reproducing formula.
Originally, these techniques  were introduced in the setting of homogeneous
spaces, and in this context they showed to be useful, for instance,
for the proof of the
$T(b)$ theorem \cite{DJS} and in the study of Trieble-Lizorkin and
Besov spaces \cite{HJTW}, \cite{HS}.

Let us denote by $\{S_k\}_{k\in\Z}$ the sequence of operators which constitute
our approximation of the identity (see Section \ref{secsk} for the precise
definition of these operators), so that for $f\in L^p(\mu)$, $1<p<\infty$,
$S_k f\to f$ in $L^p(\mu)$ as $k\to+\infty$. In this paper we will prove
estimates of the type
\begin{equation} \label{quasiorp}
\|f\|_{L^p(\mu)} \approx
\biggl\| \Bigl(\sum_{k\in\Z} |D_k f|^2 \Bigr)^{1/2}\biggr\|_{L^p(\mu)},
\end{equation}
where $D_k=S_k-S_{k-1}$, $1<p<\infty$, and the notation $A\approx B$ means
that there is some
constant $C>0$ such that $C^{-1}\,A\leq B \leq C\,A$. Notice that for $p=2$
the equation above can be rewritten as
\begin{equation} \label{quasiort}
\|f\|_{L^2(\mu)}^2 \approx \sum_{k\in\Z}\|D_k f\|_{L^2(\mu)}^2.
\end{equation}
This estimate will be a fundamental ingredient in our proof of the $T(1)$
theorem. It implies that, in some sense, the $L^2(\mu)$ decomposition
$f = \sum_{k\in\Z} D_kf$
is {\em quasiortogonal}.

\vv
In order to state the $T(1)$ theorem, we need to
introduce some notation and definitions.
Throughout all the paper we will
assume that $\mu$ is a Radon measure on $\R^d$ satisfying \rf{creix}.

\begin{definition} \label{defczo}
A kernel $k(\cdot,\cdot):\R^d\times\R^d\to\R$
is called a ($n$-dimensional) Calder\'on-Zyg\-mund (CZ) kernel if
\begin{itemize}
\item[(1)] $\ds |k(x,y)|\leq \frac{C_1}{|x-y|^n}\quad$ if $x\neq y$,
\item[(2)] there exists $0<\delta\leq1$ such that
$$|k(x,y)-k(x',y)| + |k(y,x)-k(y,x')|
\leq C_2\,\frac{|x-x'|^\delta}{|x-y|^{n+\delta}}$$
if $|x-x'|\leq |x-y|/2$.
\end{itemize}
We say that $T$ is a Calder\'on-Zygmund operator (CZO) associated to the kernel
$k(x,y)$ if for any compactly supported function $f\in L^2(\mu)$
\begin{equation} \label{eq***}
Tf(x) = \int k(x,y)\,f(y)\,d\mu(y)\qquad \mbox{if $x\not\in\supp(\mu)$}.
\end{equation}
\end{definition}

The integral in \rf{eq***} may be non convergent for $x\in\supp(\mu)$, even
for ``very nice'' functions, such as $\CC^\infty$ functions with compact
support. For this reason it is convenient to introduce the truncated operators
$T_\ve$, $\ve>0$:
$$T_\varepsilon f(x) = \int_{|x-y|>\varepsilon} k(x,y)\,f(y)\,d\mu(y).$$
It is easy to see that now this integral is absolutely convergent for any
$f\in L^2(\mu)$ and $x\in\R^d$.

We say that $T$ is bounded on $L^2(\mu)$ if the operators $T_\ve$ are bounded
on $L^2(\mu)$ uniformly on $\ve>0$.

Given a fixed constant $\rho>1$, we say that $f\in L^1_{loc}(\mu)$
belongs to the space $\bmo_\rho(\mu)$ if for some constant $C_3$
$$\sup_{Q} \frac{1}{\mu(\rho Q)}
\int_Q |f-m_Q(f)|\,d\mu\leq C_3,$$
with the supremum taken over all the cubes $Q$. By a cube $Q$ we mean a closed
cube with sides parallel to the axes and centered at some point of
$\supp(\mu)$. Also, $\rho Q$ is cube concentric with $Q$ whose side length is
$\rho$ times the side length of $Q$, and $m_Q(f)$ stands for the mean of $f$
over $Q$ with respect to $\mu$, that is, $m_Q(f) = \int_Q f\,d\mu/\mu(Q)$.

\begin{definition}
We say that $T$ is weakly bounded
(or $\rho$-weakly bounded) if
\begin{equation} \label{weakbdd}
\bigl|\langle T_\ve\chi_Q,\,\chi_Q\rangle\bigr| \leq C\,\mu(\rho Q)
\end{equation}
for any cube $Q$, uniformly on $\ve>0$.
\end{definition}

For this definition we have followed \cite{NTV3}. Let us notice that it
differs slightly from the usual definition of
weak boundedness in the doubling situation. However, the definition above
seems more natural in our context. For a discussion regarding this question,
see Section 1 of \cite{NTV3}.

Now we are ready to state the $T(1)$ theorem:

\begin{theorem} \label{t1}
If $T$ is a CZO which is weakly bounded and $T_\ve(1),\,T^*_\ve(1)\in
\bmo_\rho(\mu)$ uniformly on $\ve>0$
for some $\rho>1$, then $T$ is bounded on $L^2(\mu)$.
\end{theorem}

Some remarks are in order.
In the theorem, $T^*_\ve$ stands for the adjoint of $T_\ve$ with respect to
the duality
$\langle f,\,g\rangle = \int f\,g\,d\mu$. On the other hand, $T_\ve$ and
$T^*_\ve$ can be extended to $L^\infty(\mu)$ functions in the usual way.
The arguments are only a slight variant from the
ones of the classical doubling case. See \cite[p.300]{Stein}, for example.

Let us remark that in the case of $\mu$ being the Lebesgue measure on $\R^d$,
and also in homogeneous spaces, it has been more usual to state the
$T(1)$ theorem not
in terms of the truncated operators $T_\ve$, but in terms of some abstract
extension of $T$ to the whole space $L^2(\mu)$, which it is assumed to be
bounded from $\mathcal S$ to $\mathcal S'$ a priori. Our approach to the
$T(1)$ theorem in terms of $T_\ve$'s avoids the technical difficulties
originated from the convergence of the integral in \rf{eq***} for
$x\in\supp(\mu)$.

The proof of Theorem \ref{t1} will follow quite closely the scheme of the
proof of the $T(b)$ theorem on homogeneous spaces in \cite{DJS}.
In general, the estimates will be more difficult than in the
homogeneous case, because of the poor regularity of the measure $\mu$.
We will apply the methods developed in \cite{Tolsa3} and \cite{Tolsa5}.
In particular, the space $\rbmo(\mu)$ introduced in \cite{Tolsa3} will play
a fundamental role in the proof.

\vv
The plan of the paper is the following.
In Section \ref{lipa} we sketch the arguments for obtaining
Littlewood-Paley type estimates with respect to $\mu$.
Sections \ref{seclattice}, \ref{sec4} and \ref{secsk} deal
with the geometric construction that is needed to implement this
Littlewood-Paley theory. In Section \ref{sec6} we apply this construction to
obtain estimates such as \rf{quasiorp} and \rf{quasiort}. The rest of the paper
is devoted to the proof of the $T(1)$ theorem. First, a technical lemma
corresponding to the case $T(1)=T^*(1)=0$ is
proved in Section \ref{sect10}, and finally the theorem in its complete form is
obtained in Section \ref{sec8}, by means of the construction of a suitable
paraproduct.

\vv
Throughout all the paper the letter $C$ will be used for constants that may
change from one occurrence
to another. Constants with subscripts, such as $C_1$, do not change in
different occurrences.


\section{A Calder\'on type reproducing formula}  \label{lipa}

In this section we will describe the construction based on Coifman's ideas
that will allow the
introduction of Littlewood-Paley techniques in $L^2(\mu)$ for a measure $\mu$
satisfying \rf{creix} and non doubling in general.

We will consider a sequence of integral operators $\{S_k\}_{k\in\Z}$ given by
kernels $s_k(x,y)$ defined on $\R^d\times\R^d$. This sequence of operators
will give some kind of approximation of the identity, with $S_k \to I$ as
$k\to+\infty$ and $S_k\to 0$ as $k\to-\infty$ strongly in $L^2(\mu)$
(we say that $S_k \to S$ strongly in $L^2(\mu)$ if for any $f\in L^2(\mu)$,
$S_k f\to Sf$ in $L^2(\mu)$).
For each $x$, the support of
$s_k(x,\cdot)$ will be ``near'' some cube of scale $k$ centered at $x$, and
similarly for each $y$ the support of $s_k(\cdot,y)$ will be ``near'' some
cube of scale $k$ centered at $y$ (thus $S_kf$ approximates $f$ at
some scale $k\in\Z$). Moreover, the kernels $s_k(x,y)$ will
satisfy some appropriate size and regularity conditions and
\begin{equation}  \label{intsk}
\begin{split}
& \int s_k(x,y)\, d\mu(x) = 1 \quad \mbox{for each $y\in\supp(\mu)$,}\\
& \int s_k(x,y)\, d\mu(y) = 1 \quad \mbox{for each $x\in\supp(\mu)$.}
\end{split}
\end{equation}

For each $k$ we set $D_k= S_k - S_{k-1}$, and then, at least formally,
\begin{equation}  \label{sum2}
I=\sum_{k\in\Z} D_k
\end{equation}
We will prove that
\begin{equation} \label{sum3}
C^{-1}\,\sum_k\|D_k f\|_{L^2(\mu)}^2 \leq \|f\|_{L^2(\mu)}^2
\leq C\,\sum_k\|D_k f\|_{L^2(\mu)}^2
\end{equation}
for any $f\in L^2(\mu)$.
Now we are going to sketch the arguments for proving these inequalities, always
at a formal level.

To prove the left inequality in \rf{sum3} it is enough to show that the
operator
$\sum_k D_k^* D_k$ is bounded on $L^2(\mu)$, since
$\sum_k \|D_kf\|_{L^2(\mu)}^2 = \langle \sum_k D_k^* D_k f,\, f\rangle$.

To get the right inequality in \rf{sum3} we operate as follows.
By \rf{sum2} we have
\begin{multline}  \label{sum**}
I= \Biggl( \sum_{k\in\Z} D_k\Biggr) \,\Biggl( \sum_{j\in\Z} D_j\Biggr)
= \sum_{k\in\Z} \sum_{j\in\Z} D_{k+j} \,D_j \\
= \sum_{|k|\leq N} \sum_{j\in\Z} D_{k+j} \,D_j +
\sum_{|k|> N} \sum_{j\in\Z} D_{k+j} \,D_j.
\end{multline}
We denote $E_k = \sum_{j\in\Z} D_{k+j}\,D_j$ and $\Phi_N = \sum_{|k|\leq N}
E_k$.
Observe that if we set $D_k^N = \sum_{j:\,|j-k|\leq N} D_j$, then we also have
$\Phi_N = \sum_{k\in\Z} D_k^N \,D_k.$

Notice that in \rf{sum**} we only have stated $I= \Phi_N + (I-\Phi_N)$.
We can guess
that under the appropriate conditions, $\Phi_N\to I$ as $N\to +\infty$.
We will show that indeed this convergence occurs in the operator norm
of $L^2(\mu)$. Then,
for $N$ big enough, $\|I-\Phi_N\|_{2,2} \leq 1/2$
(where $\|\cdot\|_{2,2}$ stands for the operator norm in $L^2(\mu)$)
and so $\Phi_N$ is an
invertible operator on $L^2(\mu)$. This implies $\|f\|_{L^2(\mu)} \leq C\,
\|\Phi_N f\|_{L^2(\mu)}$ for any $f\in L^2(\mu)$.

Therefore, to see that the right inequality of \rf{sum3} holds we only have to
show that $\|\Phi_N f\|_{L^2(\mu)}^2 \leq C\,\sum_k\|D_k f\|_{L^2(\mu)}^2$.
This follows by a converse H\"older inequality argument. Given $g\in L^2(\mu)$,
we have
\begin{eqnarray}  \label{sum4}
|\langle \Phi_N f,\,g\rangle | & = & \Bigl| \sum_k \langle D_k^N D_k
f,\,g\rangle \Bigr| \, =
\, \Bigl| \sum_k \langle D_k f,\,D_k^{N*} g\rangle \Bigr|  \nonumber \\
& \leq & \sum_k \|D_k f\|_{L^2(\mu)}\,\|D_k^{N*} g\|_{L^2(\mu)} \nonumber\\
& \leq & \Bigl( \sum_k\|D_k f\|_{L^2(\mu)}^2 \Bigr)^{1/2} \,
\Bigl( \sum_k\|D_k^{N*} g\|_{L^2(\mu)}^2 \Bigr)^{1/2}.
\end{eqnarray}
From the definition of $D_k^N$ and the left inequality of \rf{sum3} we obtain
\begin{equation} \label{sum5}
\sum_k \|D_k^{N*} g\|_{L^2(\mu)}^2 \leq C\,N^2 \sum_k \|D_k^* g\|_{L^2(\mu)}^2
\leq C\,N^2\|g\|_{L^2(\mu)}^2
\end{equation}
(in our construction, we will have $D_k^* = D_k$).
Thus, by \rf{sum4} and \rf{sum5} the right inequality in \rf{sum3} follows.

\vspace{2mm}
One of the difficulties for implementing the arguments above when $\mu$ is a
non doubling measure arises from the non trivial construction of the kernels
$s_k(x,y)$ satisfying the required properties.

In case $\mu$ is doubling and satisfies $\mu(B(x,r))\approx r^n$ for all $x\in
\supp(\mu)$ and all $r>0$, the argument used by David, Journ\'e and Semmes
\cite{DJS} for homogeneous spaces works:
we fix a smooth radial function $\vphi:\R^d\lra\R$
such that $\chi_{B(0,1)} \leq \vphi \leq \chi_{B(0,2)}$ and then for each
$y\in\supp(\mu)$ and $k\in\Z$ we set
$$\vphi_{y,k}(x) = \frac{1}{r^n} \,\vphi\left(\frac{y-x}{r}\right),$$
with $r=2^{-k}$. We consider the kernel $\wt{s}_k(x,y) = \vphi_{y,k}(x)$
and so we have
\begin{equation*}
\begin{split}
& \int \wt{s}_k(x,y)\, d\mu(x) \approx 1 \quad \mbox{for each $y\in\supp(\mu)$,}\\
& \int\wt{s}_k(x,y)\, d\mu(y) \approx 1 \quad \mbox{for each $x\in\supp(\mu)$.}
\end{split}
\end{equation*}
In the estimates for proving \rf{sum3} it is essential that
$$\int d_k(x,y)\,d\mu(x) = \int d_k(x,y)\,d\mu(y) = 0.$$
So we cannot simply take $s_k(x,y):= \wt{s}_k(x,y)$.

The solution of \cite{DJS} is the following. Let $\wt{S}_k$ be the integral
operator with kernel $\wt{s}_k(x,y)$, $M_k$ the operator of multiplication
by $1/\wt{S}_k1$, and $W_k$ the operator of multiplication by $\left[
\wt{S}_k^*(1/\wt{S}_k1)\right]^{-1}$. Then we set
$S_k = M_k \wt{S}_k W_k \wt{S}_k^* M_k$. Thus the kernel of $S_k$ is
$$s_k(x,y) = \int \frac{1}{\wt{S}_k1(x)}\, \wt{s}_k(x,z)\,
\left[\wt{S}_k^*(1/\wt{S}_k1)(z)\right]^{-1}\,\wt{s}_k(y,z)\,
\frac{1}{\wt{S}_k1(y)}\, d\mu(z).$$
It is easily seen that $S_k1=1$ and, since $s_k(x,y) = s_k(y,x)$, both
identities in \rf{intsk} are satisfied.

When $\mu$ is a non doubling measure we will follow a similar approach. The
difficult step consists of obtaining functions $\vphi_{y,k}(x)$ such that
\begin{equation}  \label{eq11}
\begin{split}
& \int \vphi_{y,k}(x)\, d\mu(x) \approx 1 \quad \mbox{for each
$y\in\supp(\mu)$,}\\
& \int \vphi_{y,k}(x)\, d\mu(y) \approx 1 \quad \mbox{for each
$x\in\supp(\mu)$.}
\end{split}
\end{equation}
Nevertheless, in \cite{Tolsa5} some functions fulfilling \rf{eq11} have
been constructed. A variant of the arguments of \cite{Tolsa5} will yield the
required functions $\vphi_{y,k}$. Then we will apply the arguments of
\cite{DJS}:
we will set $S_k= M_k\wt{S}_k W_k \wt{S}_k^* M_k$, with the same notations as
above. Let us remark that, unlike in the preceding case of $\mu$ doubling,
now we will have $\vphi_{x,k}(y)\neq \vphi_{y,k}(x)$ in general, and so
$\wt{S}_k\neq \wt{S}_k^*$.

The rest of the argument for proving \rf{sum3} (which will show that all
the manipulations above dealing with the operators $D_k$ are correct) is
based on estimates analogous to the ones of \cite{DJS}, although in general
they will be more involved. One has to keep in mind that our ``dyadic'' cubes
of the $k$th scale, $k\in\Z$, will not be cubes of side length $2^{-k}$.
In the  ``dyadic'' lattice that we will construct there will not be a
direct relation between the scale $k$ of some cube $Q$ and $\mu(Q)$ or
the side length of $Q$, $\ell(Q)$.


\section{The lattice of cubes}  \label{seclattice}


\subsection{Preliminaries}   \label{secprelim}

We will assume that the constant $C_0$ in \rf{creix} has been chosen
big enough so that for all the cubes $Q\subset \R^d$ we have
$\mu(Q) \leq C_0\, \ell(Q)^n$.

\begin{definition}
Given $\alpha>1$ and $\beta>\alpha^n$, we say that the cube $Q\subset \R^d$
is $(\alpha,\beta)$-doubling if $\mu(\alpha Q) \leq \beta\,\mu(Q)$.
\end{definition}

\begin{rem} \label{moltdob}
As shown in \cite{Tolsa3},
due to the fact that $\mu$ satisfies the growth condition \rf{creix}, there
are a lot ``big'' doubling cubes. To be precise, given any point
$x\in\supp(\mu)$ and $c>0$, there exists some
$(\alpha,\beta)$-doubling cube $Q$ centered at $x$ with $l(Q)\geq c$. This
follows easily from \rf{creix} and the fact that $\beta>\alpha^n$.

On the other hand, if $\beta>\alpha^d$, then for $\mu$-a.e. $x\in\R^d$ there
exists a sequence of $(\alpha,\beta)$-doubling cubes $\{Q_k\}_k$ centered
at $x$ with $\ell(Q_k)\to0$ as $k\to \infty$.
So there are a lot of ``small'' doubling cubes too.

For definiteness, if $\alpha$ and $\beta$ are not specified, by a doubling
cube we mean a $(2,2^{d+1})$-doubling cube.
\end{rem}

Given cubes $Q,R\subset \R^d$, we denote by $z_Q$ the center of $Q$, and by
$Q_R$ the smallest cube
concentric with $Q$ containing $Q$ and $R$.

\begin{definition} \label{deltaqr}
Given two cubes $Q,R\subset\R^d$, we set
$$\delta(Q,R) = \max \left( \int_{Q_R\setminus Q}
\frac{1}{|x-z_Q|^n}\, d\mu(x),
\,\,\int_{R_Q\setminus R} \frac{1}{|x-z_R|^n}\, d\mu(x) \right).$$
\end{definition}

Notice that $\ell(Q_R) \approx \ell(R_Q) \approx \ell(Q) + \ell(R) +
\dist(Q,R)$, and if $Q\subset R$, then $R_Q=R$ and $\ell(R)\leq \ell(Q_R) \leq
2\ell(R)$.

We may treat points $x\in\supp(\mu)$ as if they were cubes (with
$\ell(x)=0$). So for $x,y\in\supp(\mu)$ and some cube $Q$, the notations
$\delta(x,Q)$ and $\delta(x,y)$ make sense. In some way, they are particular
cases of Definition \ref{deltaqr}. Of course, it may happen
$\delta(x,Q) = \infty$ or $\delta(x,y)=\infty$.

The coefficients $\delta(Q,R)$ have already appeared in our previous works
\cite{Tolsa3} and \cite{Tolsa5}. In particular, the definition of the space
$\rbmo(\mu)$ in \cite{Tolsa3} is given in terms of these coefficients:

\begin{definition} \label{defrbmo}
We say that some function $f\in L^1_{loc}(\mu)$
belongs to the space $\rbmo(\mu)$ if there exists some constant $C_4$
such that for any {\em doubling} cube
$$\frac{1}{\mu(Q)}
\int_Q |f-m_Q(f)|\,d\mu\leq C_4,$$
and for any two {\em doubling} cubes $Q\subset R$,
$$|m_Q f- m_R f|\leq C_4\,(1+\delta(Q,R)).$$
The minimal constant $C_4$ equals the $\rbmo(\mu)$ norm of $f$, which we will
denote by $\|f\|_*$.
\end{definition}

In the following lemma, proved in \cite{Tolsa5}, we recall some useful
properties of $\delta(\cdot,\cdot)$.

\begin{lemma}  \label{propdelta}
The following properties hold:
\begin{itemize}
\item[(a)] If $\ell(Q) \approx \ell(R)$ and $\dist(Q,R)\lesssim \ell(Q)$, then
$\delta(Q,R) \leq C$. In particular, $\delta(Q,\rho Q)\leq C_0\,2^n\,\rho^n$ for
$\rho>1$.

\item[(b)] Let $Q\subset R$ be concentric cubes such that there are no doubling
cubes of the form $2^k Q$, $k\geq0$, with $Q\subset 2^k Q\subset R$. Then,
$\delta(Q,R)\leq C_5$.

\item[(c)] If $Q\subset R$, then
$$\delta(Q,R) \leq C\,\left(1+\log\frac{\ell(R)}{\ell(Q)}\right).$$

\item[(d)] If $P\subset Q \subset R$, then
$$\bigl|\delta(P,R) - [ \delta(P,Q) + \delta(Q,R) ] \bigr| \leq \ve_0.$$
That is, with a different notation,
$\delta(P,R) = \delta(P,Q) + \delta(Q,R) \pm \ve_0$.
If $P$ and $Q$ are concentric, then
$\ve_0=0$: $\delta(P,R) = \delta(P,Q) + \delta(Q,R)$.

\item[(e)] For $P,Q,R\subset \R^d$,
$$\delta(P,R) \leq C_6\, + \delta(P,Q) + \delta(Q,R).$$
\end{itemize}
\end{lemma}

The constants that appear in (b), (c), (d) and (e) depend on $C_0,n,d$.
The constant $C$ in (a) depends, further, on the constants that are
implicit in the relations $\approx,\, \lesssim$.

Let us insist on the fact that a notation such as $a= b \pm \ve$ does not mean
any precise equality but the estimate $|a-b|\leq \ve$.

Notice that if we set $D(Q,R) = 1 + \delta(Q,R)$ for $Q\neq R$ and
$D(Q,Q) = 0$, then $D(\cdot,\cdot)$ is a quasidistance on the set of cubes,
by (e) in the preceding lemma.

If we denote by $\wt{Q}$ the smallest doubling cube of the form $2^kQ$,
$k\geq0$, by (b) we know that $\wt{Q}$ is not far from $Q$ (using the
quasidistance $D$). So $Q$ and $\wt{Q}$ may have very different sizes, but we
still have $D(Q,\wt{Q})\leq C$.

\vv
In Remark \ref{moltdob} we have explained that there a lot of big and small
doubling cubes. In the following lemma we state a more precise result about the
existence of small doubling cubes in terms of $\delta(\cdot,\cdot)$.

\begin{lemma}\label{precisio}
There exists some (big) constant $\gamma_0>0$ depending only on $C_0$, $n$
and $d$
such that if $R_0$ is some cube centered at some point of $\supp(\mu)$ and
$\alpha>\gamma_0$, then
for each $x\in R_0\cap \supp(\mu)$ such that $\delta(x,2R_0) >\alpha$ there
exists some doubling cube $Q\subset 2R_0$ centered at $x$ satisfying
\begin{equation}  \label{err1}
|\delta(Q,2R_0) - \alpha|\leq \ve_1,
\end{equation}
where $\ve_1$ depends only on $C_0$, $n$ and $d$ (but not on $\alpha$).
\end{lemma}

See \cite{Tolsa5} again for the proof of this lemma.

As in (d) of Lemma \ref{propdelta},
instead of \rf{err1}, often we will write $\delta(Q,2R_0) = \alpha
\pm\ve_1$.

Now we are going to state a similar result concerning the existence of big
doubling cubes with some precise estimate involving the ``distance''
$\delta(\cdot,\cdot)$.

\begin{lemma}\label{precisio'}
There exists some (big) constant $\gamma_0>0$ depending only on $C_0$, $n$
and $d$ such that for any fixed $\alpha>\gamma_0$,
if $R_0$ is some cube centered at some point of $\supp(\mu)$ with
$\delta(R_0,\R^d)>\alpha$, then there
exists some doubling cube $S \supset R_0$ concentric with $R_0$, with
$\ell(S)\geq 2\ell(R_0)$, satisfying
\begin{equation}  \label{err1'}
|\delta(R_0,S) - \alpha|\leq \ve_1,
\end{equation}
where $\ve_1$ depends only on $C_0$, $n$ and $d$ (but not on $\alpha$).
\end{lemma}

The proof follows by arguments analogous to the ones for proving
Lemma \ref{precisio}.

For convenience, we will assume that the constant $\ve_1$ of Lemmas
\ref{precisio} and \ref{precisio'} has been chosen so that $\ve_1\geq\ve_0$.


\subsection{Cubes of different generations}  \label{sscubdif}

\begin{definition}  \label{defstop}
We say that $x\in \supp(\mu)$ is a {\em stopping point (or stopping cube)} if
$\delta(x,Q) <\infty$ for some cube $Q\ni x$ with $0<\ell(Q)<\infty$.
We say that $\R^d$ is a {\em initial cube} if $\delta(Q,\R^d)<\infty$
for some cube $Q\ni x$ with $0<\ell(Q)<\infty$.
The cubes $Q$ with $0<\ell(Q)<\infty$ are called {\em transit cubes}.
\end{definition}

It is easily seen that if $\delta(x,Q) <\infty$ for some transit cube $Q$
containing $x$, then $\delta(x,Q') <\infty$ for any other transit cube $Q'$
containing $x$. Also, if $\delta(Q,\R^d) <\infty$ for
some transit cube $Q$, then $\delta(Q',\R^d) <\infty$ for any transit cube
$Q'$.

Notice that the points (which are also cubes following our convention)
which are not
stopping cubes have not received any special name. The same happens for
$\R^d$ if it is not an initial cube. This is because this cubes will not
play any specific role in our geometric construction.

We will take some big positive integer $A$
whose precise value will be fixed after knowing or choosing several additional
constants. In particular, we assume that $A$ is much bigger than the constants
$\ve_0,\,\ve_1$ and $\gamma_0$ of Section \ref{secprelim}.

Now we are ready to introduce the definition of generations of cubes (in a
first case).

\begin{definition}  \label{defcubs}
Assume that $\R^d$ is not an initial cube. We fix some doubling cube
$R_0\subset \R^d$. This will be our ``reference'' cube.
For each integer $j\geq 1$ we let $R_{-j}$ be some doubling cube concentric
with $R_0$, containing $R_0$, and such that $\delta(R_0,R_{-j})=jA \pm \ve_1$
(which exists because of Lemma \ref{precisio'}).
If $Q$ is a transit cube,
we say that $Q$ is a cube of generation $k\in \Z$ if it is a doubling cube and
for some cube $R_{-j}$ containing $Q$ we have
\begin{equation}  \label{identk}
\delta(Q,R_{-j}) = (j + k)A \pm \ve_1.
\end{equation}
If $Q\equiv x$ is a stopping cube, we say that $Q$ is a cube of generation $k$
if for some cube $R_{-j}$ containing $x$ we have
\begin{equation} \label{ineqk}
\delta(x,R_{-j}) \leq (j + k)A \pm \ve_1.
\end{equation}
\end{definition}

Notice that the cubes $R_{-j}$, $j\geq1$, are cubes of generation $-j$
and that if $Q$ is a transit cube of generation $k$ contained in {\em some}
$R_{-j}$, then $\delta(Q,R_{-j}) =
(j+k)A \pm 3\ve_1$ (with ``$\leq$'' if $Q$ is a stopping cube), by (d) of Lemma
\ref{propdelta}. So, in some way, modulo some small errors, the chosen
reference $R_{-j}$ does not matter.

Observe that if $\R^d$ is not an initial cube, then for any $x\in\supp(\mu)$
there are cubes of all generations
$k\in \Z$ centered at $x$. Indeed, for $A$ big enough we have
$\ell(R_{-j})\to+\infty$ as $j\to+\infty$.
So for any $x\in \supp(\mu)$ we choose
$R_{-j}$ such that $x\in \frac{1}{2}R_{-j}$,
and then we only have to apply Lemma \ref{precisio}.

For any $x\in\supp(\mu)$,
we denote by $Q_{x,k}$ some fixed doubling cube centered $x$
of the $k$th generation. If $x$ is not a stopping point
and $\R^d$ is not an initial cube, then all the cubes will be
transit cubes and the identity \rf{identk} holds for them.
If $x$ is a stopping point, then there exists some $k_x\in\Z$ such that all the
cubes of generations $k<k_x$ centered at $x$ are transit cubes,
and all the cubes centered at $x$ of generation $k>k_x$ coincide with the
point $x$ (we can think they have ``collapsed'' in the point $x$).

In case $\R^d$ is an initial cube we have to modify a little the definition
above because not all the cubes $R_{-j}$ in that definition exist.

\begin{definition}  \label{defcubs'}
Assume that $\R^d$ is an initial cube. Then we choose $\R^d$ as our
``reference'': If $Q$ is a transit cube,
we say that $Q$ is a cube of generation $k\geq1$ if
\begin{equation}  \label{identk'}
\delta(Q,\R^d) = k A\pm \ve_1.
\end{equation}
If $Q\equiv x$ is a stopping cube, we say that $Q$ is a cube of generation
$k\geq1$ if
\begin{equation} \label{ineqk'}
\delta(x,\R^d) \leq k A \pm \ve_1.
\end{equation}
Moreover, for all $k\leq 1$ we say that $\R^d$ is a cube of generation $k$.
\end{definition}

As in the case where $\R^d$ is not an initial cube, for any $x$ we also have
cubes of all generations centered at $x$ (we have to think that $\R^d$ is
centered at all the points $x\in\supp(\mu)$).

Observe that the last definition coincides with Definition \ref{defcubs} with
the convention (that we will follow) $R_{-j}=\R^d$ for $j\geq0$.

\begin{definition}
For any $x\in\supp(\mu)$,
we denote by $Q_{x,k}$ some fixed cube centered $x$
of the $k$th generation.
\end{definition}

If $Q_{x,k}\neq\{x\}$ and $Q_{x,k}\subset R_{-j}$,
then we have $\delta(Q_{x,k},R_{-j})\approx (j+k)A$, because
$A$ is much bigger than $\ve_1$.
However, the estimate \rf{identk} is much sharper.
This will very useful in our construction.

The constants $\ve_0$ and $\ve_1$ should be understood as upper
bounds for some ``errors'' and deviations of our construction from the
classical dyadic lattice.

\vv
It is easily seen that if $A$ is big enough, then  $\ell(Q_{x,k+1})\leq
\ell(Q_{x,k})/10$. So $\ell(Q_{x,k}) \to 0$ as $k\to +\infty$.
In fact, the following more precise result holds.

\begin{lemma}  \label{mides2} \label{mides}
If we take $A$ is big enough, then there exists some $\eta>0$ such that
if $x,y\in\supp(\mu)$ are such that $2Q_{x,k} \cap 2Q_{y,k+m} \neq \varnothing$
(with $m\geq1$),
then $\ell(Q_{y,k+m}) \leq 2^{-\eta m}\,\ell(Q_{x,k})$.
\end{lemma}

See \cite{Tolsa5} for the proof.


\section{The construction of the functions $\vphi_{y,k}$}  \label{sec4}

In this section we will explain how to construct the functions
$\vphi_{y,k}$ which will
originate the kernels $s_k(x,y)$. This construction will follow the same
lines as the one in \cite{Tolsa5}, although with some simplifications.

We denote
$$\sigma:= 100 \ve_0 + 100\ve_1 + 12^{n+1} C_0.$$

We introduce two new constants $\alpha_1, \alpha_2>0$ whose precise value will
be fixed below. For the moment, let us say that $\ve_0, \ve_1, C_0,
\ll \sigma \ll \alpha_1\ll\alpha_2\ll A$.

\begin{definition} Let $y\in \supp(\mu)$. If $Q_{y,k}$ is a transit cube,
we denote by $Q_{y,k}^1$, $\wh{Q}_{y,k}^1$, $Q_{y,k}^2$, $\wh{Q}_{y,k}^2$,
$Q_{y,k}^3$ some
doubling cubes centered at $y$ containing $Q_{y,k}$ such that
\begin{equation}  \label{igu}
\begin{split}
& \delta(Q_{y,k},Q_{y,k}^1) = \alpha_1 \pm \ve_1,\\
& \delta(Q_{y,k},\wh{Q}_{y,k}^1) =  \alpha_1 +\sigma \pm \ve_1,\\
& \delta(Q_{y,k},Q_{y,k}^2) = \alpha_1 + \alpha_2 \pm \ve_1,\\
& \delta(Q_{y,k},\wh{Q}_{y,k}^2) = \alpha_1 + \alpha_2 +\sigma \pm \ve_1 ,\\
& \delta(Q_{y,k},Q_{y,k}^3) = \alpha_1 + \alpha_2 +2\,\sigma \pm \ve_1.
\end{split}
\end{equation}
By Lemma \ref{precisio} and the definitions of Subsection \ref{sscubdif},
we know that all these cubes exist.

If $Q_{y,k}=\R^d$, we set $Q^1_{y,k}= \wh{Q}_{y,k}^1= Q_{y,k}^2 =
\wh{Q}_{y,k}^2 = Q_{y,k}^3 =\R^d$. If $Q_{y,k}\equiv y$ is a stopping cube and
$Q_{y,k-1}\equiv y$ is also a stopping cube, we set
$Q^1_{y,k}= \wh{Q}_{y,k}^1= Q_{y,k}^2 =
\wh{Q}_{y,k}^2 = Q_{y,k}^3 =y$. If $Q_{y,k}\equiv y$ is a stopping cube but
$Q_{y,k-1}\equiv y$ is not, then we choose $Q_{y,k}^1$, $\wh{Q}_{y,k}^1$, $Q_{y,k}^2$, $\wh{Q}_{y,k}^2$,
$Q_{y,k}^3$ so that they are contained in $Q_{y,k-1}$, centered at $y$ and
\begin{equation}  \label{igu'}
\begin{split}
& \delta(Q_{y,k}^1,Q_{y,k-1}) = A- \alpha_1 \pm \ve_1,\\
& \delta(\wh{Q}_{y,k}^1,Q_{y,k-1}) = A- \alpha_1 -\sigma \pm \ve_1,\\
& \delta(Q_{y,k}^2,Q_{y,k-1}) = A- \alpha_1 - \alpha_2 \pm \ve_1,\\
& \delta(\wh{Q}_{y,k}^2,Q_{y,k-1}) = A- \alpha_1 - \alpha_2 -\sigma \pm \ve_1 ,\\
& \delta(Q_{y,k}^3,Q_{y,k-1}) = A -\alpha_1 - \alpha_2 -2\,\sigma \pm \ve_1.
\end{split}
\end{equation}
If any of these cubes does not exists because $\delta(y,Q_{y,k-1})$ is not
big enough, we let this cube be the point $\{y\}$.
\end{definition}

If $Q_{y,k}$ is a transit cube, then the identities \rf{igu'} are also
satisfied by $Q_{y,k}^1$, $\wh{Q}_{y,k}^1$, $Q_{y,k}^2$, $\wh{Q}_{y,k}^2$,
$Q_{y,k}^3$, by (d) in Lemma \ref{propdelta}.
So in this case it would be possible to define $Q_{y,k}^1$,
$\wh{Q}_{y,k}^1$, $Q_{y,k}^2$, $\wh{Q}_{y,k}^2$,
$Q_{y,k}^3$ by the identities \rf{igu} too. However, we think that
the definition is more clear if we take $Q_{y,k}$ as the reference, as in
\rf{igu}.

\begin{lemma}  \label{construc}
Let $y\in \supp(\mu)$. If we choose the constants
$\alpha_1$, $\alpha_2$ and  $A$ big enough, we have
\begin{equation}  \label{inclu}
Q_{y,k}\subset Q_{y,k}^1 \subset \wh{Q}_{y,k}^1 \subset Q_{y,k}^2
\subset \wh{Q}_{y,k}^2\subset Q_{y,k}^3 \subset Q_{y,k-1}.
\end{equation}
\end{lemma}

The proof of this lemma follows from an easy calculation. See \cite{Tolsa5}
for the details.

For a fixed $k$, cubes of the $k$th generation
may have very different sizes for different $y$'s.
The same happens for the cubes $Q_{y,k}^1$ and $Q_{y,k}^2$. Nevertheless, in
\cite{Tolsa5} it has been shown that we still have some kind
of regularity:

\begin{lemma}  \label{regu}
Given $x,y\in\supp(\mu)$,
let $Q_x,\,Q_y$ be cubes centered at $x$ and $y$ respectively, and assume that
$Q_x\cap Q_y\neq \varnothing$ and that there exists some cube $R_0$
containing $Q_x\cup Q_y$, with
$|\delta(Q_x,R_0) - \delta(Q_y,R_0)| \leq 10\ve_1$. If $R_y$ is some cube
centered at $y$ containing $Q_y$ with
$\delta(Q_y,R_y) \geq \sigma - 10\ve_1$, then $Q_x\subset R_y$.
As a consequence, we have:
\begin{enumerate}
\item[(a)] If $Q_{x,k}^1\cap Q_{y,k}^1\neq\varnothing$, then
$Q_{x,k}^1\subset \wh{Q}_{y,k}^1$,
in particular $x\in \wh{Q}_{y,k}^1$.
\item[(b)] If $Q_{x,k}^2\cap Q_{y,k}^2\neq\varnothing$,then $Q_{x,k}^2\subset \wh{Q}_{y,k}^2$,
in particular $x\in \wh{Q}_{y,k}^2$.
\item[(c)] If $Q_{x,k} \cap Q_{y,k}\neq \varnothing$, then
$Q_{x,k} \subset Q_{y,k-1}$.
\end{enumerate}
\end{lemma}

So, although we cannot expect to have the equivalence
$$y\in Q_{x,k}^1 \Leftrightarrow x\in Q_{y,k}^1,$$
we still have something quite close to it, because
the cubes $Q_{x,k}^1$ and $\wh{Q}_{x,k}^1$
are close one each other in the quasimetric $D(\cdot,\cdot)$, since
$\delta(Q_{x,k}^1,\wh{Q}_{x,k}^1)$ is small (at least in front of $A$).
Of course, the same idea applies if we change $1$ by $2$ in the superscripts
of the cubes.

Now we are going to define the functions $\vphi_{y,k}$. First we introduce
the auxiliary functions $\psi_{y,k}$.

\begin{definition}  \label{defpsi}
For any $y\in \supp(\mu)$,
the function $\psi_{y,k}$ is a function such that
\begin{enumerate}
\item $0\leq \psi_{y,k}(x)\leq \min\left(\dfrac{4}{\ell(Q_{y,k}^1)^n},\,
\dfrac{1}{|y-x|^n}\right)$,
\item $\psi_{y,k}(x) = \dfrac{1}{|x-y|^n}$ if $x\in \wh{Q}_{y,k}^2\setminus
Q_{y,k}^1$,
\item $\supp(\psi_{y,k})\subset Q_{y,k}^3$,
\item $|\psi_{y,k}'(x)|\leq C_7\,
\min\left(\dfrac{1}{\ell(Q_{y,k}^1)^{n+1}},\,  \dfrac{1}{|y-x|^{n+1}}\right)$.
\end{enumerate}
\end{definition}

It is not difficult to check that such a function exists if we  choose $C_7$
big enough.
We have to take into account that $2\wh{Q}_{y,k}^2 \subset Q_{y,k}^3$.
This is due to the fact that $\delta(\wh{Q}_{y,k}^2,2\wh{Q}_{y,k}^2)\leq 4^nC_0
< \delta(\wh{Q}_{y,k}^2,Q_{y,k}^3)$ if $\ell(\wh{Q}_{y,k}^2)\neq0$.

In the definition of $\psi_{y,k}$, if $Q_{y,k}^1=\{y\}$, then one must take
$1/\ell(Q_{y,k}^1)=\infty$. If $\wh{Q}_{y,k}^2 = \{y\}$, then we set
$\psi_{y,k}\equiv0$. If $Q_{y,k}=\R^d$, we set $\psi_{y,k}\equiv0$.
These choices satisfy the conditions in the definition
of $\psi_{y,k}$ stated above.

\begin{definition}  \label{defphi}
For all $y\in\supp(\mu)$, we set
$\vphi_{y,k}(x) = \alpha_2^{-1}\, \psi_{y,k}(x)$.
\end{definition}

Choosing $\alpha_2$ big enough, the largest part of the $L^1(\mu)$ norm of
$\psi_{y,k}$ and $\vphi_{y,k}$ will come from the integral over
$Q_{y,k}^2\setminus \wh{Q}_{y,k}^1$.
We state this in a precise way in the following lemma.

\begin{lemma}  \label{propor}
There exists some constant $\ve_2$ depending on $n$, $d$, $C_0$, $\ve_0$,
$\ve_1$ and $\sigma$ (but not on $\alpha_1$, $\alpha_2$ nor $A$) such that
if $Q_{y,k}^1 \neq \{y\},\,\R^d$, then
\begin{equation}  \label{err30}
\left| \|\psi_{y,k}\|_{L^1(\mu)} - \alpha_2 \right| \leq \ve_2
\end{equation}
and
\begin{equation}  \label{err31}
\left| \|\psi_{y,k}\|_{L^1(\mu)} - \int_{Q_{y,k}^2\setminus \wh{Q}_{y,k}^1}
 \frac{1}{|y-x|^n}\,d\mu(x) \right| \leq \ve_2.
\end{equation}
\end{lemma}

The proof of this result is an easy calculation that we will skip.
A direct consequence of it is
$$\lim_{\alpha_2\to \infty} \frac{1}{\alpha_2}\int_{Q_{y,k}^2\setminus
\wh{Q}_{y,k}^1}
\frac{1}{|y-x|^n}\,d\mu(x) = 1$$
for $y\in\supp(\mu)$ such that $Q_{y,k}^1 \neq \{y\},\R^d$.

In order to study some of the properties of the functions
$\vphi_{y,k}$, we need to introduce some additional notation.

\begin{definition}
Let $x\in\supp(\mu)$ and assume that $Q_{x,k}\neq \R^d$.
We denote by $\QH_{x,k}$ a doubling cube centered at
$x$ such that $\delta(\QH_{x,k},Q_{x,k-1})= A - \alpha_1- \alpha_2 - 3\,\sigma
\pm \ve_1$.
Also, we denote by $\check{Q}_{x,k}^1$ and $\Check{\Check{Q}}_{x,k}^1$ some
doubling cubes centered at $x$ such that
\begin{equation*}
\begin{split}
& \delta(\check{Q}_{x,k}^1,Q_{x,k-1}) = A - \alpha_1 + \sigma \pm \ve_1,\\
& \delta(\Check{\Check{Q}}_{x,k}^1,Q_{x,k-1}) = A - \alpha_1 + 2\sigma \pm \ve_1
\end{split}
\end{equation*}
(the idea is that the symbols $\,\,\wh{ }\,\,$ and $ \,\,\check{ }\,\,$
are inverse operations, modulo some small errors).
If any of the cubes
$\check{Q}_{x,k}^1,\Check{\Check{Q}}_{x,k}^1, \QH_{x,k}$
does not exist because $\delta(x,Q_{x,k-1})$ is not big enough, then we let
it be the point $x$. If $Q_{x,k}=\R^d$, then we set
$\QH_{x,k} = \check{Q}_{x,k}^1 = \Check{\Check{Q}}_{x,k}^1 = \R^d$.
\end{definition}

So when $Q_{x,k}$ is a transit cube, we have
\begin{equation*}
\begin{split}
&\delta(\QH_{x,k},Q_{x,k})= \alpha_1+ \alpha_2 + 3\,\sigma \pm \ve_1,\\
&\delta(\check{Q}_{x,k}^1,Q_{x,k}) = \alpha_1 - \sigma \pm \ve_1,\\
&\delta(\Check{\Check{Q}}_{x,k}^1,Q_{x,k}) = \alpha_1 - 2\sigma \pm \ve_1.
\end{split}
\end{equation*}
and one should think that $\QH_{x,k}$ is
a cube a little bigger than
$\wh{Q}^3_{x,k}$, while $\check{Q}_{x,k}^1$ is a little smaller than
$Q_{x,k}^1$. Also, $\Check{\Check{Q}}_{x,k}^1$ is a little smaller than
$\check{Q}_{x,k}^1$, but still much bigger than $Q_{x,k}$.

\begin{lemma} \label{phi}
Let $x,y\in \supp(\mu)$. For $\alpha_1$ and $\alpha_2$ big enough, we have:
\begin{itemize}
\item[(a)] If $x\in Q_{x_0,k}$ and $y\not\in \wh{Q}_{x_0,k}^3$, then
$\vphi_{y,k}(x) = 0$. In particular,
$\vphi_{y,k}(x) = 0$ if $y\not\in \wh{Q}_{x,k}^3$.

\item[(b)] If $y\in \check{Q}_{x,k}^1$, then $\ds \vphi_{y,k}(x) \leq C\,
\frac{\alpha_2^{-1}}{\ell(\check{Q}_{x,k}^1)^n}.$

\item[(c)] For all $y\in \R^d$,
$$\vphi_{y,k}(x) \leq \frac{\alpha_2^{-1}}{|y-x|^n},$$
and if $y\in Q_{x,k}^2 \setminus \wh{Q}_{x,k}^1$, then
$$\ds \vphi_{y,k}(x) = \frac{\alpha_2^{-1}}{|y-x|^{n}}.$$

\item[(d)] If $x\in Q_{x_0,k}$, then $$|\ds \vphi_{y,k}'(x)| \leq
C\, \alpha_2^{-1} \,
\min\left(\frac{1}{\ell(\check{Q}_{x_0,k}^1)^{n+1}}, \, \frac{1}{|y-x|^{n+1}}
\right).$$
\end{itemize}
\end{lemma}

Notice that, in Definition \ref{defpsi} of the functions $\psi_{y,k}$,
the properties that define these functions are stated with respect to cubes
centered at $y$ ($Q_{y,k}^1$, $Q_{y,k}^2$, $Q_{y,k}^3$...). In this lemma
some analogous properties are stated, but these properties have to do with
cubes centered at $x$ or containing $x$
($Q_{x_0,k}$, $\check{Q}_{x,k}^1$, $Q_{x,k}^2$, $\wh{Q}_{x,k}^3$...).

\begin{proof}
\begin{itemize}

\item[(a)]
Let $x_0\in \supp(\mu)$ and $x\in Q_{x_0,k}$.
If $\vphi_{y,k}(x)\neq0$, then
$x\in Q_{y,k}^3$. Then $Q_{x_0,k}^3 \cap Q_{y,k}^3\neq\varnothing$
and so
$y\in Q_{y,k}^3\subset \wh{Q}_{x_0,k}^3$ (as in Lemma \ref{regu}).

\vv
\item[(b)] Let $y\in \check{Q}_{x,k}^1$. We know that
$$\vphi_{y,k}(x) \leq C\,\alpha_2^{-1}\,
\frac{1}{\ell(Q_{y,k}^1)^{n}}.$$
So we are done if we see that $\ell(Q_{y,k}^1) \geq \ell(\check{Q}_{x,k}^1)$.

As in Lemma \ref{regu}, we have
$$y\in \check{Q}_{x,k}^1 \Rightarrow \check{Q}_{y,k}^1 \cap\check{Q}_{x,k}^1
\neq \varnothing
\Rightarrow \check{Q}_{x,k}^1\subset Q_{y,k}^1.$$
Thus $\ell(\check{Q}_{x,k}^1)\leq \ell(Q_{y,k}^1).$

\vv
\item[(c)] The first inequality follows from the definition of $\psi_{y,k}$
and $\vphi_{y,k}$. The second statement is also straightforward. Indeed,
if $y\in Q_{x,k}^2\setminus \wh{Q}_{x,k}^1$, then by Lemma \ref{regu} we get
$x\in \wh{Q}_{y,k}^2\setminus Q_{y,k}^1.$ Notice that, in particular, this
implies $\wh{Q}_{y,k}^2\neq \{y\},\R^d$. We only
have to look at the definitions of $\psi_{y,k}$ and $\vphi_{y,k}$ again.

\vv
\item[(d)]
Suppose that $y\in \check{Q}_{x_0,k}^1$. In this case we must show that
$$|\vphi_{y,k}'(x)| \leq
C\, \frac{\alpha_2^{-1}}{\ell(\check{Q}_{x_0,k}^1)^{n+1}}.$$
It is enough to see that
$\ell(Q_{y,k}^1) \geq \ell(\check{Q}_{x_0,k}^1)$. This follows from the inclusion
$Q_{y,k}^1 \supset \check{Q}_{x_0,k}^1$, which holds because
$y\in \check{Q}_{y,k}^1 \cap \check{Q}_{x_0,k}^1$ and then we can apply Lemma
\ref{regu}.

On the other hand, since by definition we have
$$|\vphi_{y,k}'(x)| \leq C\,\frac{\alpha_2^{-1}}{|y-x|^{n+1}},$$
we are done.
\end{itemize}
\end{proof}

Some of the estimates in the preceding lemma will be used to prove
next result, which was one of our main goals in this section.

\begin{lemma} \label{convo}
For any $\ve_3>0$, if $\alpha_1$ and $\alpha_2$ are big
enough, for all $z_0\in \supp(\mu)$ we have
\begin{equation} \label{convo1}
\int \vphi_{z_0,k}(x)\,d\mu(x)\leq 1 + \ve_3 \quad \mbox{ and }\quad
\int \vphi_{y,k}(z_0)\,d\mu(y) \leq 1 + \ve_3.
\end{equation}
If $z_0\in\supp(\mu)$ is such that there exists some transit cube
$Q_k$ of the $k$th generation with $Q_k\ni z_0$, then
\begin{equation} \label{convo2}
1-\ve_3 \leq \int \vphi_{z_0,k}(x)\,d\mu(x)\quad \mbox{ and }\quad
1-\ve_3\leq \int \vphi_{y,k}(z_0)\,d\mu(y).
\end{equation}
\end{lemma}

\begin{proof}
Let us see  \rf{convo2} first. So we assume that there exist
some transit cube $Q_{k}$ of the $k$th generation containing $z_0$.
Since $z_0\in Q_{k}\subset \check{Q}_{k}^1$, we have
$\check{Q}_{k}^1 \subset Q_{z_0,k}^1$. In particular, $\ell(Q_{z_0,k}^1)>0$.
So the inequality
$$1-\ve_3 \leq \int \vphi_{z_0,k}(x)\,d\mu(x)$$
is a direct consequence of Lemma \ref{propor}.

We consider now the second inequality in \rf{convo2}.
By Lemma \ref{propor} and the second equality of (c) in Lemma \ref{phi}
we get
\begin{eqnarray*}
\int \vphi_{y,k}(z_0)\,d\mu(y) & \geq & \int_{Q_{z_0,k}^2\setminus
\wh{Q}_{z_0,k}^1}
 \vphi_{y,k}(z_0)\,d\mu(y)\\
& \geq & \int_{Q_{z_0,k}^2\setminus \wh{Q}_{z_0,k}^1}
\frac{\alpha_2^{-1}}{|y-z_0|^n} \,d\mu(y)\\
& \geq & \alpha_2^{-1}\,(\alpha_2-2\ve_2).
\end{eqnarray*}
So the second inequality in \rf{convo2} holds if we take $\alpha_2$ big enough.

Consider now \rf{convo1}. The first estimate follows easily from the
definitions
\ref{defpsi} and \ref{defphi}. Let us see the second inequality of \rf{convo1}.
By (a) in Lemma \ref{phi} have
$$\int \vphi_{y,k}(z_0)\,d\mu(y) = \int_{\wh{Q}_{z_0,k}^3}
\vphi_{y,k}(z_0)\,d\mu(y).$$
Thus we can write
\begin{equation}  \label{desc}
\int \vphi_{y,k}(z_0)\,d\mu(y) = \int_{\wh{Q}_{z_0,k}^3\setminus
\check{Q}_{z_0,k}^1}
\vphi_{y,k}(z_0)\,d\mu(y) + \int_{\check{Q}_{z_0,k}^1} \vphi_{y,k}(z_0)
\,d\mu(y).
\end{equation}

Let us estimate the first integral on the right hand side of \rf{desc}.
Using the first inequality in (c) of Lemma \ref{phi} we obtain
\begin{eqnarray}  \label{desc3}
\int_{\wh{Q}_{z_0,k}^3\setminus\check{Q}_{z_0,k}^1} \vphi_{y,k}(z_0)\,d\mu(y)
& \leq & \int_{\wh{Q}_{z_0,k}^3\setminus\check{Q}_{z_0,k}^1}
\frac{\alpha_2^{-1}}{|y-z_0|^n}\, d\mu(y) \nonumber \\
& = & \delta(\check{Q}_{z_0,k}^1,\wh{Q}_{z_0,k}^3) \,\alpha_2^{-1}
\nonumber\\
& \leq & \alpha_2^{-1}\,
(\alpha_2 + 4\,\sigma + 2\,\ve_1).
\end{eqnarray}

Let us consider the last integral in \rf{desc}
(only in the case $\check{Q}_{z_0,k}^1\neq\{z_0\},\R^d$).
By (b) in Lemma \ref{phi} we have
\begin{equation} \label{desc4}
\int_{\check{Q}_{z_0,k}^1} \vphi_{y,k}(z_0)\,d\mu(y)
\leq \int_{\check{Q}_{z_0,k}^1}
\frac{C\,\alpha_2^{-1}}{\ell(\check{Q}_{z_0,k}^1)^n}\,d\mu(y)
\leq C\,C_0\,\alpha_2^{-1}.
\end{equation}
From \rf{desc3} and \rf{desc4} we get \rf{convo1}.
\end{proof}


\section{The kernels $s_k(x,y)$} \label{secsk}

In this section we will introduce the operators $S_k$ mentioned in Section
\ref{lipa} and we will obtain some estimates for their kernels $s_k(x,y)$.

We will assume
that we have chosen $\ve_3=1/2$ in Lemma \ref{convo}. Recall that then
$1/2\leq \int\vphi_{y_0,k}(x) d\mu(x)\leq 3/2$ and
$1/2\leq \int\vphi_{y,k}(x_0) d\mu(y) \leq3/2$ if $Q_{x_0,k}$ and
$Q_{y_0,k}$ are transit cubes.

\begin{definition}
Let $f\in L^1_{loc}(\mu)$ and $x\in\supp(\mu)$. If $Q_{x,k}\neq \R^d$,
then we set
$$\wt{S}_k f(x) = \int \vphi_{y,k}(x)\,f(y)\,d\mu(y) +
\max\biggl( 0,\,\frac{1}{4}- \int \vphi_{y,k}(x)\,d\mu(y)\biggr)\, f(x).$$
\end{definition}

Observe that, formally, $\wt{S}_k$ is an integral operator with the following
positive kernel:
\begin{equation} \label{kerlio1}
\wt{s}_k(x,y) = \vphi_{y,k}(x) +
\max\biggl( 0,\,\frac{1}{4}- \int \vphi_{y,k}(x)\,d\mu(y)\biggr)
\,\delta_{x}(y),
\end{equation}
where $\delta_x$ is the Dirac delta at $x$.
If $Q_{x,k}$ is a transit cube, by Lemma \ref{convo} we have
$$\wt{S}_k f(x) = \int \vphi_{y,k}(x)\,f(y)\,d\mu(y).$$
Notice also that for {\em all} $x\in\supp(\mu)$ we have
$1/4\leq \wt{S}_k1(x)\leq 3/2$.

Now we can define the operators $S_k$:

\begin{definition}
Assume that $Q_{x,k}\neq \R^d$ for some $x\in \supp(\mu)$.
Let $M_k$ be the operator of multiplication by
$m_k(x):= 1/\wt{S}_k1(x)$ and $W_k$ the operator of multiplication by
$w_k(x) := 1/ \wt{S}_k^*(1/\wt{S}_k1)$.
We set $S_k := M_k\,\wt{S}_k\,W_k\,\wt{S}_k^*\,M_k$.
If $Q_{x,k}=\R^d$ for some $x\in \supp(\mu)$, then we set $S_k:=0$.
\end{definition}

Observe that if $Q_{x,k}$ and $Q_{y,k}$ are transit cubes, then
$$S_k f(x) = \int s_k(x,y)\,f(y)\,d\mu(y),$$
where $s_k(\cdot,\cdot)$ is the kernel
\begin{equation}  \label{sskk}
s_k(x,y) = \int m_k(x)\, \wt{s}_k(x,z)\,w_k(z)\,\wt{s}_k(y,z)\,m_k(y)\, d\mu(z).
\end{equation}

The following estimates are a direct consequence of the
statements in Lemma \ref{convo} and the definitions above.

\begin{lemma} \label{transla}
For all $k\in\Z$, if $x\in\supp(\mu)$ is such that
$Q_{x,k}\neq\R^d$, then $2/3\leq m_k(x) \leq 4$ and $0\leq w_k(x)\leq 6$.
\end{lemma}

\begin{proof}
As mentioned above, $1/4\leq\wt{S}_k1(x)\leq 3/2$ and so
$2/3\leq m_k(x) \leq 4$.

On the other hand, we also have $\wt{S}_k^*1(x)\geq1/4$,
and then
$$\wt{S}^*_k(1/\wt{S}_k 1)(x) \geq \frac{2}{3} \,\wt{S}^*_k(1) \geq
\frac{1}{6},$$
and so $w_k(x) \leq 6$.
\end{proof}

In the following lemma we show the localization, size and regularity
properties that fulfil the kernels $s_k(x,y)$.

\begin{lemma} \label{losire}
For each $k\in\Z$ the following properties hold:
\begin{itemize}
\item[(a)] If $Q_{x,k}$ is a transit cube,
then $\supp(s_k(x,\cdot)) \subset Q_{x,k-1}$.

\item[(b)] If $Q_{x,k}$ and $Q_{y,k}$ are transit cubes, then
\begin{equation} \label{sum30}
0\leq s_k(x,y) \leq \frac{C}{(\ell(Q_{x,k}) + \ell(Q_{y,k}) + |x-y|)^n}.
\end{equation}

\item[(c)] If $Q_{x,k}$, $Q_{x',k}$, $Q_{y,k}$ are transit cubes, and
$x,x'\in Q_{x_0,k}$ for some $x_0\in\supp(\mu)$, then
\begin{equation} \label{sum33}
|s_k(x,y) - s_k(x',y)| \leq C\, \frac{|x-x'|}{\ell(Q_{x_0,k})} \cdot
\frac{1}{(\ell(Q_{x,k}) + \ell(Q_{y,k}) + |x-y|)^n}.
\end{equation}
\end{itemize}
\end{lemma}

\begin{proof}
{\bf (a)} From \rf{sskk} we see that if $s_k(x,y)\neq 0$, then there exists
some $z\in\supp\mu$ such that $\vphi_{z,k}(x)\neq0$ and
$\vphi_{z,k}(y)\neq0$. Thus $z\in \wh{Q}^3_{x,k}\cap \wh{Q}^3_{y,k}$, and from
Lemma \ref{regu} we get $y\in\QH_{x,k}\subset Q_{x,k-1}$.

\vv
\noi {\bf (b)}
By \rf{sskk} and Lemma \ref{transla} we have
$$s_k(x,y) \leq C\,\int \wt{s}_k(x,z)\,\wt{s}_k(y,z)\,d\mu(z).$$
Since $\wt{s}_k(x,z)=\vphi_{z,k}(x)\leq C/\ell(Q_{x,k})^n$, we get
$$s_k(x,y) \leq \frac{C}{\ell(Q_{x,k})^n} \int \vphi_{z,k}(y)\,d\mu(y)
\leq \frac{C}{\ell(Q_{x,k})^n}.$$
Similarly it can shown that $s_k(x,y) \leq  C/\ell(Q_{y,k})^n$. So it only
remains to see that $s_k(x,y) \leq C/|x-y|^n$.

Recall that $\wt{s}_k(x,z)\leq C/|x-z|^n$ and $\wt{s}_k(y,z)\leq C/|y-z|^n$.
Then we have
\begin{eqnarray*}
s_k(x,y)  & \leq &
C\,\int_{|x-z|\geq |x-y|/2} \wt{s}_k(x,z)\,\wt{s}_k(y,z)\,d\mu(z)\\
&&\mbox{} + C\,\int_{|x-z|<|x-y|/2} \wt{s}_k(x,z)\,\wt{s}_k(y,z)\,d\mu(z) \\
&\leq & \frac{C}{|x-y|^n} \int\wt{s}_k(y,z)\,d\mu(z)  +
\frac{C}{|x-y|^n} \int\wt{s}_k(x,z)\,d\mu(z)\\
&\leq & \frac{C}{|x-y|^n}.
\end{eqnarray*}

\vv
\noi {\bf (c)}
Using Lemma \ref{transla} we get
\begin{eqnarray*}
|s_k(x,y) - s_k(x',y)| &\leq& C\, |m_k(x) - m_k(x')| \int
\wt{s}_k(x,z)\,\wt{s}_k(y,z)\,d\mu(z) \\
&&\mbox{} + C \int |\wt{s}_k(x,z) - \wt{s}_k(x',z)|\, \wt{s}_k(y,z) d\mu(z)\\
&=& A + B.
 \end{eqnarray*}
Let us estimate the term $A$. Operating as in (b), we obtain
$$\int
\wt{s}_k(x,z)\,\wt{s}_k(y,z)\,d\mu(z) \leq
\frac{C}{(\ell(Q_{x,k}) + \ell(Q_{y,k}) + |x-y|)^n}.$$
On the other hand, since $x,x'\in Q_{x_0,k}$, $\wt{S}_k 1\approx 1$, and
\begin{equation}\label{sum32}
|\vphi_{z,k}'(w)|\leq \frac{C}{(\ell(\check{Q}^1_{x_0,k}) + |w-z|)^{n+1}}
\end{equation}
for all $w\in Q_{x_0,k}$, we get
\begin{eqnarray*}
|m_k(x)- m_k(x')| & \leq  & \left|\int C\,(\vphi_{z,k}(x) -
\vphi_{z,k}(x'))\,d\mu(z)\right| \\
& \leq & \int \frac{C\,|x-x'|}{(\ell(Q_{x_0,k}) + |x-z|)^{n+1}}\,d\mu(z)
\, \leq \, C\,\frac{|x-x'|}{\ell(Q_{x_0,k})}.
\end{eqnarray*}
So $A$ verifies inequality \rf{sum33}.

Let us consider the term $B$ now. By \rf{sum32} we obtain
\begin{eqnarray*}
B & \leq & \int \frac{C\,|x-x'|}{(\ell(\check{Q}^1_{x_0,k}) + |x-z|)^{n+1}}\,
\wt{s}_k(y,z)\,d\mu(z) \\
& = & \int_{|z-y|\geq|x-y|/2} + \int_{|z-y|<|x-y|/2}  \, =\,B_1 + B_2
\end{eqnarray*}
Since $\vphi_{z,k}(y) \leq C/(\ell(\check{Q}^1_{y,k}) + |y-z|)^n$, we have
\begin{eqnarray*}
B_1 & \leq & \frac{C\,|x-x'|}{(\ell(\check{Q}^1_{y,k}) + |x-y|)^n}
\int \frac{1}{(\ell(Q_{x_0,k}) + |x-z|)^{n+1}}\,d\mu(z) \\
& \leq & C\, \frac{|x-x'|}{\ell(Q_{x_0,k})}\cdot
\frac{1}{(\ell(\check{Q}^1_{y,k}) + |x-y|)^n}.
\end{eqnarray*}
It is easy to check that
$\ell(Q_{x,k})\leq 2(\ell(\check{Q}^1_{y,k}) + |x-y|)$. Indeed if $|x-y| \leq
\ell(Q_{x,k})/2$, then $y\in Q_{x,k}$ and so $Q_{x,k}\subset\check{Q}^1_{y,k}$
and so $\ell(Q_{x,k})\leq \ell(\check{Q}^1_{y,k})$. Thus the term $B_1$ also
satisfies \rf{sum33}.

Let us turn our attention to $B_2$. In this case we have
\begin{eqnarray*}
B_2 & \leq & C\,\frac{|x-x'|}{(\ell(\check{Q}^1_{x_0,k}) + |x-y|)^{n+1}}\,
\int \wt{s}_k(y,z)\,d\mu(z) \\
& \leq & C\, \frac{|x-x'|}{\ell(Q_{x_0,k})}\cdot
\frac{1}{(\ell(\check{Q}^1_{x_0,k}) + |x-y|)^n}.
\end{eqnarray*}
Thus we only have to check that $\ell(Q_{x,k}) + \ell(Q_{y,k}) \leq
C\,(\ell(\check{Q}^1_{x_0,k}) + |x-y|)$. Because $x\in Q_{x_0,k}$, we have
$Q_{x,k}\subset \check{Q}^1_{x_0,k}$ and so $\ell(Q_{x,k}) \leq
\ell(\check{Q}^1_{x_0,k})$. Let us see that $\ell(Q_{y,k}) \leq C\,
(\ell(\check{Q}^1_{x_0,k}) + |x-y|)$. If $|x_0-y|\geq \ell(Q_{y,k})/2$, then
$$\frac{1}{2}\, \ell(Q_{y,k}) \leq |x-x_0| + |x-y| \leq
C\,\ell(Q_{x_0,k}) + |x-y|.$$
If $|x_0-y|<\ell(Q_{y,k})/2$, then $x_0\in Q_{y,k}$ and so $Q_{y,k}\subset
\check{Q}_{x_0,k}^1$, which yields $\ell(Q_{y,k}) \leq
\ell(\check{Q}^1_{x_0,k})$.
\end{proof}

Notice that, in general, the functions $\wt{s}_k(x,y) = \vphi_{y,k}(x)$
do not have any smoothness with respect to the variable $y$. However, the
kernels $s_k(x,y)$ defined above have regularity in both variables, because
$s_k(y,x) = s_k(x,y)$. On the other hand, this smoothness appears to be
somewhat weaker than the regularity in $x$ of the functions $\vphi_{y,k}(x)$.

\begin{rem}  \label{kerlio}
Taking the (formal) definition \rf{kerlio} of the kernels
$\wt{s}_k(x,y)$, it is easily seen that
the properties of the kernels $s_k(x,y)$ in (a), (b) and (c) of the
lemma above also hold without the assumptions without assuming that
$Q_{x,k}$, $Q_{x',k}$ and $Q_{y,k}$ are transit cubes.
Indeed, the statements are trivial is any of these cubes coincides with $\R^d$,
and if any of them is a stopping cube, then it is not difficult to
check that all the estimates in the proof above are also valid.
\end{rem}


\section{Littlewood-Paley type estimates}  \label{sec6}

We recall some notation introduced in Section \ref{lipa}.
For each $k\in\Z$, we set $D_k = S_k - S_{k-1}$,
$E_k = \sum_{j\in\Z} D_{k+j}\,D_j$ and, for each $N\geq1$,
$\Phi_N = \sum_{|k|\leq N} E_k$.

Notice that $D_k1=0$ for all $k\in\Z$ except in the case $k=1$ for $\R^d$
being an initial cube.

\begin{lemma} \label{lemadk}
We have:
\begin{itemize}
\item[(a)] $\|D_j\,D_k\|_{2,2} \leq C\,2^{-|j-k|\,\eta}$ for all $j,k\in\Z$ and
some $\eta>0$.
\item[(b)] $\sum_{k\in Z} D_k = I$, with strong convergence in $L^2(\mu)$.
\item[(c)] The series $\sum_{j\in\Z} D_{k+j}\,D_j =:E_k$ converges strongly
in $L^2(\mu)$ and $$\|E_k\|_{2,2} \leq C\,|k|\,2^{-|k|\eta}$$ for all
$k\in\Z$.
\item[(d)] $\Phi_N\to I$ as $N\to+\infty$ in the operator norm in $L^2(\mu)$.
\end{itemize}
\end{lemma}

\begin{proof}
For simplicity we assume that all the cubes $Q_{x,k}$,
$x\in\supp(\mu),\,k\in\Z$, are transit cubes.
In the final part of the proof we will give some hints for the general case.
Moreover, we only have to prove the assertion (a). The others follow from (a)
by the Cotlar-Knapp-Stein Lemma, as in \cite{DJS}.

Assume $j\geq k+2$. The kernel of the operator $D_j\,D_k$ is given by
$$K_{j,k}(x,y) = \int d_j(x,z)\,d_k(z,y)\,d\mu(z).$$
Since $\supp(d_j(x,\cdot))\subset Q_{x,j-2}$, we have
$$|K_{j,k}(x,y)| \leq \int_{z\in Q_{x,j-2}}
|d_j(x,z)\,(d_k(z,y) - d_k(x,y))|\,d\mu(z).$$
By (b) of Lemma \ref{losire} (taking into account that $Q_{x,j-2}\subset
Q_{x,k}$),
$$|d_k(z,y) - d_k(x,y)| \leq C\,\frac{\ell(Q_{x,j-2})}{\ell(Q_{x,k})} \cdot
\frac{1}{(\ell(Q_{x,k}) + \ell(Q_{y,k}) + |x-y|)^n}.$$
By Lemma \ref{mides2} we have $\ell(Q_{x,j-2}) \leq
C\,2^{-\eta|j-k|}\,\ell(Q_{x,k})$ for some $\eta>0$. Therefore,
\begin{eqnarray}  \label{eqxxx}
|K_{j,k}(x,y)| & \leq & C\,2^{-\eta|j-k|}\,
\frac{1}{(\ell(Q_{x,k}) + \ell(Q_{y,k}) + |x-y|)^n}\int |d_j(x,z)|\,d\mu(z)
\nonumber \\
& \leq &  C\,2^{-\eta|j-k|}\,
\frac{1}{(\ell(Q_{x,k}) + \ell(Q_{y,k}) + |x-y|)^n}.
\end{eqnarray}
Also, we have $\supp(K_{j,k}(x,\cdot)) \subset Q_{x,k-3}$ and
$\supp(K_{j,k}(\cdot,y)) \subset Q_{y,k-3}$. Indeed, if $K_{j,k}(x,y)\neq 0$
then there exists some $z\in Q_{x,j-2}\cap Q_{y,k-2}$, and so
$y\in Q_{x,k-3}$ and $x\in Q_{y,k-3}$.
Then we obtain
\begin{eqnarray}\label{eqwww21}
\int |K_{j,k}(x,y)|\,d\mu(y) & \leq & C\,2^{-\eta|j-k|}\, \int_{Q_{x,k-3}}
\frac{1}{(\ell(Q_{x,k}) + |x-y|)^n}\,d\mu(y) \nonumber \\
& \leq & C\,2^{-\eta|j-k|}\,(1+\delta(Q_{x,k},Q_{x,k-3})) \leq
C\,2^{-\eta|j-k|}.
\end{eqnarray}
In an analogous way, we get
\begin{equation}\label{eqwww22}
\int |K_{j,k}(x,y)|\,d\mu(x) \leq C\,2^{-\eta|j-k|}.
\end{equation}
Therefore, by Schur's Lemma we have $\|D_j\,D_k\|_{p,p} \leq C\,2^{-\eta|j-k|}$
for all $p\in[1,\infty]$ if $j\geq k+2$.

On the other hand, for $k\geq j+2$, operating in a similar way, we also obtain
$\|D_j\,D_k\|_{p,p} \leq C\,2^{-\eta|j-k|}$,
and if $|j-k|\leq1$, then we have $\|D_j\,D_k\|_{p,p} \leq \|D_j\|_{p,p} \,
\|D_k\|_{p,p}\leq C$. Thus the assertion (a) of the lemma holds in any case.

\vv
If there exist stopping cubes, then by Remark \ref{kerlio} the kernels of
the operators $S_k$ satisfy properties which are similar to the ones stated in
Lemma \ref{losire}, and some estimates as the ones above work.
If $\R^d$ is an initial cube, then $\int d_1(x,y)\,d\mu(y) \neq 0 $, in general.
However in the arguments above it is used $\int d_j(x,y)\,d\mu(y) = 0 $ only to
estimate $\|D_j\,D_k\|_{p,p}$ in the case $j\geq k+2$, and notice that $D_k=0$ for
$k\leq 0$.
\end{proof}

By the estimates of the preceding lemma and by a new application of
Cotlar-Knapp-Stein Lemma, arguing as in Section \ref{lipa},
we get:

\begin{theorem}  \label{teolp}
If $f\in L^2(\mu)$, then
$$C^{-1}\,\sum_k\|D_k f\|_{L^2(\mu)}^2 \leq \|f\|_{L^2(\mu)}^2
\leq C\,\sum_k\|D_k f\|_{L^2(\mu)}^2.$$
\end{theorem}

We omit the detailed proof of this result. We only have to apply the same
arguments as in \cite{DJS} (see also \cite{HJTW}).
From this theorem we derive the following corollaries.

\begin{coro}
Let $1<p<\infty$.
If $f\in L^p(\mu)$, then
\begin{equation} \label{lppp}
C^{-1}\, \biggl\| \Bigl(\sum_k |D_k f|^2 \Bigr)^{1/2}\biggr\|_{L^p(\mu)}
\leq \|f\|_{L^p(\mu)} \leq
C\, \biggl\| \Bigl(\sum_k |D_k f|^2 \Bigr)^{1/2}\biggr\|_{L^p(\mu)}.
\end{equation}
\end{coro}

\begin{proof}
The right inequality follows from the left one (with $p'$ instead of $p$).
Indeed, by an argument similar to the one used for $p=2$ in \rf{sum4}, it
follows that
$$\|\Phi_N f\|_{L^p(\mu)} \leq
C\, \biggl\| \Bigl(\sum_k |D_k f|^2 \Bigr)^{1/2}\biggr\|_{L^p(\mu)}.$$
In Lemma \ref{inverbmo} below we will show that $\Phi_N$ is bounded and
invertible in $L^p(\mu)$, and so $\|f\|_{L^p(\mu)} \leq
C\,\|\Phi_N f\|_{L^p(\mu)}$.

The left inequality in \rf{lppp} will be proved using techniques of vector
valued
Calder\'on-Zygmund operators. These techniques, which are standard in the
classical doubling case, have been extended by Garc\'{\i}a-Cuerva and Martell
\cite{GM} to the case of non homogeneous spaces.

Let us denote by $L^p(\ell^2,\mu)$ the Banach space of sequences of functions
$\{g_k\}_{k\in\Z}$, $g_k\in L^1_{loc}(\mu)$, such that
$$\left(\int \Bigl( \sum_k |g_k|^2\Bigr)^{p/2}\,
d\mu\right)^{1/p}< \infty.$$
Let us consider the operator $D:L^p(\mu)\!\to\! L^p(\ell^2,\mu)$ given by
$Df = \{D_k f\}_{k\in\Z}$.
By Theorem \ref{teolp}, $D$ is bounded from $L^2(\mu)$ into
$L^2(\ell^2,\mu)$.
From the results in \cite{GM}, it follows that if the kernel
$d(x,y):=\{d_k(x,y)\}_k$ of $D$ satisfies
\begin{itemize}
\item[(1)] $\ds \|d(x,y)\|_{\ell^2} \leq \frac{C}{|x-y|^n}$ for $x\neq y$, and
\item[(2)] $\ds \int_{|x-y|\geq 2|x-x'|}
\bigl( \|d(x,y) - d(x',y)\|_{\ell^2} + \|d(y,x) - d(y,x')\|_{\ell^2}\bigr)\,
d\mu(y) \leq C$,
\end{itemize}
then $D$ is bounded from $L^p(\mu)$ into $L^p(\ell^2,\mu)$, $1<p<\infty$,
because $D$ is
a {\em vector-valued Calder\'on-Zygmund operator}.
Thus we only have to check that these conditions are satisfied.

Let us see that the first one holds. Given $x,y\in\supp(\mu)$, $x\neq y$,
let $j\in\Z$ be such that $y\in Q_{x,j}\setminus Q_{x,j+1}$. Since
$\supp\, d_k(x,\cdot)\subset Q_{x,k-2}$, we have
\begin{eqnarray*}
\sum_{k\in\Z} |d_k(x,y)|^2 & \leq & \sum_{k\leq j+2} \frac{C}{(\ell(Q_{x,k}) +
|x-y|)^{2n}} \\
& \leq &  \frac{C}{|x-y|^{2n}} + \sum_{k\leq j} \frac{C}{\ell(Q_{x,k})^{2n}}
\, \leq \,\frac{C}{|x-y|^{2n}}.
\end{eqnarray*}

Now we will show that condition (2) is also satisfied. Since $d(x,y)=d(y,x)$,
we only have to deal with the term $\|d(x,y) - d(x',y)\|_{\ell^2}$.
Let $h\in\Z$ be such that $x'\in Q_{x,h} \setminus Q_{x,h+1}$, and suppose that
$y\in Q_{x,j}\setminus Q_{x,j+1}$ for some $j\leq h-10$.
Notice that $d_k(x,y) - d_k(x',y)=0$ if $k>j+4$. Indeed,
we have
$$\supp(d_k(x,\cdot) - d_k(x',\cdot))\subset Q_{x,k-2} \cup Q_{x',k-2}.$$
If  $k\geq h$, then $Q_{x,k-2} \cup Q_{x',k-2}\subset
Q_{x,h-3} \subset Q_{x,j+1}$,
and if $j+4<k<h$, then we have
$Q_{x,k-2} \cup Q_{x',k-2}\subset Q_{x,k-3} \subset Q_{x,j+1}$.

Assuming $k\leq j+4$, we get
$$|d_k(x,y)-d_k(x',y)| \leq C\, \frac{|x-x'|}{\ell(Q_{x,k})}\cdot
\frac{1}{(\ell(Q_{x,k}) + |x-y|)^n},$$
since $x'\in Q_{x,h}$, with $h>k$.
Therefore,
\begin{eqnarray} \label{lpp23}
\sum_{k\in\Z} |d_k(x,y)-d_k(x',y)|^2 & \leq &C \sum_{k\leq j+4}
\left( \frac{|x-x'|}{\ell(Q_{x,k})}\cdot
\frac{1}{(\ell(Q_{x,k}) + |x-y|)^n} \right)^2 \nonumber \\
& \leq & C\, \left(\frac{|x-x'|}{\ell(Q_{x,j+4})\, |x-y|^{n}}\right)^2.
\end{eqnarray}
Then, using condition (1) and \rf{lpp23} we obtain
\begin{multline*}
\int_{|x-y|\geq 2|x-x'|} \|d(x,y) - d(x',y)\|_{\ell^2}\, d\mu(y)\\
\begin{split}
& \leq \int_{Q_{x,h-10}\setminus B(x,2|x-x'|)} \frac{C}{|x-y|^n}\, d\mu(y)\\
&\quad + C \sum_{i=10}^\infty \int_{Q_{x,h-i-1} \setminus Q_{x,h-i}}
\frac{|x-x'|}{\ell(Q_{x,h-i+4})\, |x-y|^{n}}\, d\mu(y) \\
& \leq  C + C \sum_{i=10}^\infty\frac{|x-x'|}{\ell(Q_{x,h-i+4})} \leq C.
\end{split}
\end{multline*}
\end{proof}

\begin{coro} \label{corolprbmo}
If $f\in\rbmo(\mu)$ and $Q_k$ is a cube of generation $k\in\Z$, then
\begin{equation}  \label{lprbmo}
\sum_{j=k}^{+\infty} \|D_j f\|_{L^2(\mu|Q_k)}^2 \leq C\,\|f\|_*^2\,\mu(Q_k).
\end{equation}
\end{coro}

Recall that, by definition, we assume that the cubes $Q_k$ are doubling.

\begin{proof}
For $N$ big enough and $j\geq k+N$,
$D_j f (x) = D_j (\chi_{3/2Q_k}f) (x)$ if $x\in Q_k$. Thus
\begin{eqnarray*}
\sum_{j=k+N}^{+\infty} \|D_j f\|_{L^2(\mu|Q_k)}^2 & = &
\sum_{j=k+N}^{+\infty}
\|D_j ((f - m_{Q_k}f)\,\chi_{3/2Q_k})\|_{L^2(\mu|Q_k)}^2 \\
& \leq & C\,\|f - m_{Q_k}f\|_{L^2(\mu|3/2Q_k)}^2 \\
& \leq & C\,\|f\|_*^2\,\mu(2Q_k) \leq C\,\|f\|_*^2\,\mu(Q_k).
\end{eqnarray*}
Now we only have to check that
\begin{equation}  \label{sum49}
\|D_j f\|_{L^2(\mu|Q_k)}^2 \leq C\,\|f\|_*^2\,\mu(Q_k)
\end{equation}
for $j=k,\ldots,k+N-1$. We set
$$\|D_j f\|_{L^2(\mu|Q_k)} \leq \|S_j (f-m_{Q_k}f)\|_{L^2(\mu|Q_k)} +
\|S_{j-1} (f-m_{Q_k}f)\|_{L^2(\mu|Q_k)}.$$
For each $j$, we denote by $N_j$ the least integer such
that $Q_{x,j-1}\subset 2^{N_j}Q_{x,j}$. We have
\begin{eqnarray*}
|S_jf(x) - m_{Q_{x,j}}f| & \leq &
\left|\int s_j(x,y)\,(f(y) - m_{Q_{x,j}}f)\,d\mu(y)\right| \\
& \leq & C\,\sum_{m=1}^{N_j} \int_{2^mQ_{x,j}\setminus2^{m-1}Q_{x,j}}
\frac{|f(y)-m_{Q_{x,j}}  f|}{\ell(2^mQ_{x,j})^n}\,d\mu(y) \\
& \leq & C\,\sum_{m=1}^{N_0}\frac{1}{\ell(2^mQ_{x,j})^n}
\int_{2^mQ_{x,j}} |f(y) - m_{Q_{x,j}}f|\,d\mu(y).
\end{eqnarray*}
Since $\delta(Q_{x,j},2^mQ_{x,j})\leq C$, we get
$$\int_{2^mQ_{x,j}}  |f(y) - m_{Q_{x,j}}f|\,d\mu(y)
\leq C\,\mu(2^{m+1}Q_{x,j}),$$
and so
\begin{eqnarray*}
|S_jf(x) - m_{Q_{x,j}}f| & \leq & C\,\sum_{m=1}^{N_j}
\frac{\mu(2^{m+1}Q_{x,j})}{\ell(2^mQ_{x,j})^n}\,\|f\|_* \\
& \leq & C\,(1+\delta(Q_{x,j},Q_{x,j-1}))\,\|f\|_* \leq C\,\|f\|_*.
\end{eqnarray*}

For $j=k,\ldots,k+N-1$, since $\delta(Q_{x,j}, Q_k)\leq C\,N$, we have
$|m_{Q_k} f - m_{Q_{x,j}}f|\leq C\,\|f\|_*$. Thus
$$|S_jf(x) - m_{Q_k}f| \leq C\|f\|_*,$$
and then \rf{sum49} holds.
\end{proof}

Observe that the same arguments above show that if $f\in\rbmo(\mu)$, then
\begin{equation}  \label{lprbmo2}
\sum_{j=k-N_0}^{+\infty} \|D_j f\|_{L^2(\mu|Q_k)}^2 \leq C\,\|f\|_*^2\,
\mu(Q_k),
\end{equation}
where $N_0>0$ is some fixed integer, and $C$ depends on $N_0$ now.


\section{The $T(1)$ theorem in the case $T(1)=T^*(1)=0$}
\label{sect10}


\subsection{The main steps}
For simplicity, we will prove the $T(1)$ theorem assuming that {\em there are
no stopping cubes and $\R^d$ is not an initial cube}.
However, we claim that our arguments can be extended quite easily to the
general situation.

The kernels of the truncated operators $T_\ve$ do not satisfy the gradient
condition in the definition of CZO's. For this reason we need to introduce
the regularized operators $\wt{T}_\ve$.
Let $\vphi$ be a radial $\CC^\infty$ function with $0\leq \vphi\leq1$,
vanishing on $B(0,1/2)$ and identically equal to $1$ on $\R^d\setminus B(0,1)$.
For each $\ve>0$, we consider the integral operator $\wt{T}_\ve$ with kernel
$\vphi((x-y)/\ve)\cdot k(x,y)$.
It is easily seen that
\begin{equation}  \label{diftve}
|T_\ve f - \wt{T}_\ve f| \leq M_\mu f,
\end{equation}
where $M_\mu$ is the centered maximal Hardy-Littlewood operator.
So $T_\ve$ is bounded on $L^2(\mu)$ uniformly on $\ve>0$ if and only if
the same holds for
$\wt{T}_\ve$. The kernel of $\wt{T}_\ve$ is $L^\infty$-bounded and
it is straightforward to check that it is a CZ kernel itself,
with constants $C_1$ and $C_2$ in Definition \ref{defczo} uniform on $\ve>0$.

The following lemma will be very useful for our arguments. It shows that
the hypotheses of weak boundedness and $T_\ve(1),
T^*_\ve(1)\in \bmo_\rho(\mu)$ can be substituted by conditions about the
$L^p(\mu)$ boundedness over characteristic functions of cubes.

\begin{lemma}  \label{caracterist}
Let $T$ be a CZO. For any fixed $\rho,\gamma>1$ and $1<p<\infty$,
the following conditions are equivalent:
\begin{itemize}
\item[(a)] $T$ is weakly bounded and $T_\ve(1)\in \bmo_\rho(\mu)$ uniformly
on $\ve>0$.
\item[(b)] $T$ is weakly bounded and $T_\ve(1)\in \rbmo(\mu)$ uniformly
on $\ve>0$.
\item[(c)] For any cube $Q\subset \R^d$,
\begin{equation} \label{carac}
\|T_\ve\chi_Q\|_{L^p(\mu)} \leq C\,\mu(\gamma Q)^{1/p}
\end{equation}
uniformly on $\ve>0$.
\end{itemize}
\end{lemma}

The proof of this result follows by arguments similar to the ones in the proof
of \cite[Theorem 8.4]{Tolsa3}.

From \rf{diftve} and \rf{carac} we infer that, under the assumptions of Theorem
\ref{t1}, $\wt{T}_\ve$ is also
bounded over characteristic functions of cubes (i.e. satisfies \rf{carac})
uniformly on $\ve>0$. Of course, the same happens for $\wt{T}_\ve^*$.

In this section we will prove the following technical version of the $T(1)$
theorem.

\begin{lemma} \label{t10}
Let $k(\cdot,\cdot)$ be a CZ kernel with
$k(\cdot,\cdot)\in L^\infty(\R^d\times\R^d)$. Let $T$ be the integral operator
\begin{equation}  \label{wwwq}
Tf(x) = \int k(x,y) f(y)\,d\mu(y),\qquad f\in L^2(\mu).
\end{equation}
Assume that for $p=2$ and, in the case $n>1$, also for $p=n/(n-1)$, we have
$$\|T\chi_Q\|_{L^p(\mu)}\leq C\,\mu(2Q)^{1/p},\qquad
\|T^*\chi_Q\|_{L^p(\mu)}\leq C\,\mu(2Q)^{1/p},$$
for any cube $Q$. If moreover $T(1)=T^*(1)=0$, then $T$ is bounded on
$L^2(\mu)$,  and $\|T\|_{2,2}$ is bounded above by some constant
independent of $\|k(\cdot,\cdot)\|_\infty$.
\end{lemma}

Notice that, as we are assuming $k(\cdot,\cdot)\in L^\infty(\R^d\times\R^d)$,
from condition (1) in the definition of CZ kernel we
derive that $k(\cdot,y),\,k(x,\cdot)\in L^2(\mu)$ uniformly on $x,y$. As
a consequence,
the integral in \rf{wwwq} is convergent. So in this case, when we say that
$T$ is bounded we are not talking about the uniform boundedness of the
truncated operators $T_\ve$, but about the operator $T$ itself.
Observe also that, in particular this lemma can be applied to $\wt{T}_\ve$,
for each $\ve>0$.

In the whole section we will assume that $T$ is an operator fulfilling
the assumptions of Lemma \ref{t10}.
\vv

For each $i\in\Z$, $x\in\supp(\mu)$, we denote $u_{x,i}(z) = s_i(x,z) -
s_{i-1}(x,z) = d_i(x,z)$.

The first step of the proof of Lemma \ref{t10} consists of estimating the term
$|\langle u_{x,j},\, Tu_{y,k}\rangle|$. As we shall see, this part of the proof
will be more involved than in \cite{DJS}, basically due to the fact that the
functions $u_{x,i}$ are much less localized in our present situation.

\begin{lemma} \label{difi}
Under the assumptions of Lemma \ref{t10},  there exists some $\nu>0$
depending on $\delta,\eta$ such that for $x,y\in\supp(\mu)$ and
$j,k\in\Z$, we have
\begin{itemize}
\item[(a)] If $2Q_{x,j-3}\cap 2Q_{y,k-3} = \varnothing$, then
$$|\langle Tu_{x,j},\, u_{y,k}\rangle| \leq C\,2^{-\nu|j-k|}
\,\frac{(\ell(Q_{x,j-2}) \wedge \ell(Q_{y,k-2}))^{\delta/2}}{
(\ell(Q_{x,j-2}) + \ell(Q_{y,k-2}) + |x-y|)^{n+\delta/2}}.
$$

\item[(b)] If $2Q_{x,j-3}\cap 2Q_{y,k-3} \neq \varnothing$, then
\begin{eqnarray*}
|\langle Tu_{x,j},\, u_{y,k}\rangle| & \leq &
\frac{C\,2^{-\nu|j-k|}}{(\ell(Q_{x,j}) + |x-y|)^n}\, \chi_{Q_{x,j-7}}(y) \\
&& \mbox{} +
\frac{C\,2^{-\nu|j-k|}}{(\ell(Q_{y,k}) + |x-y|)^n}\, \chi_{Q_{y,k-7}}(x).
\end{eqnarray*}

\end{itemize}
\end{lemma}

We defer the proof of these estimates until Subsection \ref{provadifi}.
Now we will
see how from this result Lemma \ref{t10} follows by arguments analogous
to the ones of \cite{DJS}.

For each $j,k\in\Z$ we set $T_{j,k}= D_k\,T\,D_j$. The $L^2(\mu)$ norm of
$T_{j,k}$ is easily estimated by means of Lemma \ref{difi}, as we show in
next lemma.

\begin{lemma}  \label{tjk}
The operator $T_{j,k}$ is bounded on $L^2(\mu)$ with norm $\|T_{j,k}\|_{2,2}
\leq C\,2^{-\nu|j-k|}$.
\end{lemma}

\begin{proof}
The kernel of $T_{j,k}$ is given by $t_{j,k}(x,y) = \langle Tu_{x,j},\,
u_{y,k}\rangle$. We will apply Schur's Lemma, using the estimates of the
preceding lemma, interchanging $T$ and $T^*$ when necessary.

We have
$$\int |t_{j,k}(x,y)|\, d\mu(y) = \int_{y:\,2Q_{y,k-3}\cap 2Q_{x,j-3}=
\varnothing}
+ \int_{y:\,2Q_{y,k-3}\cap 2Q_{x,j-3}\neq\varnothing} = I_1 + I_2.$$
By Lemma \ref{difi} we get
$$I_1\leq  C\,2^{-\nu|j-k|}\int_{y:\,2Q_{y,k-3}\cap 2Q_{x,j-3}=\varnothing}
\frac{\ell(Q_{x,j-2})^{\delta/2}}{|x-y|^{n+\delta/2}} \,d\mu(y).$$
If $2Q_{y,k-3}\cap 2Q_{x,j-3}=\varnothing$, then $|x-y|\geq
\ell(Q_{x,j-2})/2$. Thus
$$I_1\leq C\,2^{-\nu|j-k|}\int_{|x-y|\geq\ell(Q_{x,j-2})/2}
\frac{\ell(Q_{x,j-2})^{\delta/2}}{|x-y|^{n+\delta/2}} \,d\mu(y) \leq
C\, 2^{-\nu|j-k|}.$$

We estimate $I_2$ now. By Lemma \ref{difi} we obtain
\begin{eqnarray*}
I_2 & \leq & \int_{Q_{x,j-7}}\frac{C\,2^{-\nu|j-k|}}{(\ell(Q_{x,j})
+ |x-y|)^n}\,d\mu(y) \\
&&\mbox{} + \int_{y:\,x\in Q_{y,k-7}}
\frac{C\,2^{-\nu|j-k|}}{(\ell(Q_{y,k}) + |x-y|)^n}\,d\mu(y)
\, = \, I_{2,1} + I_{2,2}.
\end{eqnarray*}
We have
\begin{eqnarray*}
I_{2,1} & \leq & \int_{Q_{x,j}}
\frac{C\,2^{-\nu|j-k|}}{\ell(Q_{x,j})^n}\,d\mu(y) +
\int_{Q_{x,j-7}\setminus Q_{x,j}} \frac{C\,2^{-\nu|j-k|}}{|x-y|^n}\,d\mu(y) \\
& \leq & C\,2^{-\nu|j-k|} (1 + \delta(Q_{x,j},Q_{x,j-7})) \, \leq
\,C\,2^{-\nu|j-k|}.
\end{eqnarray*}
Finally we turn our attention to $I_{2,2}$. Observe that if $x\in Q_{y,k-7}$,
then $y\in Q_{x,k-8}$, and so
\begin{eqnarray*}
I_{2,2} & \leq & \int_{Q_{x,k}} \frac{C\,2^{-\nu|j-k|}}{\ell(Q_{x,k})^n}\,
d\mu(y) +
\int_{Q_{x,k-8}\setminus Q_{x,k}} \frac{C\,2^{-\nu|j-k|}}{|x-y|^n}\, d\mu(y) \\
& \leq & C\,2^{-\nu|j-k|} (1 + \delta(Q_{x,k},Q_{x,k-8}) \,
\leq \, C\,2^{-\nu|j-k|}.
\end{eqnarray*}

Thus $\int |t_{j,k}(x,y)|\,d\mu(y) \leq C\,2^{-\nu|j-k|}$. By the
symmetry of the assumptions, we also have $\int |t_{j,k}(x,y)|\,d\mu(x) \leq
C\,2^{-\nu|j-k|}$. By Schur's Lemma we get
$\|T_{j,k}\|_{2,2}\leq C\,2^{-\nu|j-k|}$, and we are done.
\end{proof}

Let $J,\,K\subset \Z$ be finite sets. We set $T_{J,K} = \sum_{j\in J}
\sum_{k\in K} D_j^N\,D_j\,T\, D_k\,D_k^N$, where $N$ is an integer such that
$\|I-\Phi_N\|_{2,2} \leq 1/2$, as explained in Section \ref{lipa}. Then we have

\begin{lemma} \label{TJK}
The operator $T_{J,K}$ is bounded on $L^2(\mu)$ with $\|T_{J,K}\|_{2,2}\leq
C\,N^2$, where $C$ does not depend on $J$ or $K$.
\end{lemma}

\begin{proof}
For $f,g\in L^2(\mu)$, by Lemma \ref{tjk}, we have
\begin{eqnarray*}
|\langle T_{J,K}f,\,g\rangle| & \leq & \sum_{j\in J} \sum_{k\in K}
|\langle D_j\,T\,D_k\,D_k^N f,\, D_j^N g\rangle| \\
& \leq &
\sum_{j,k} 2^{-\nu|j-k|}\, \|D^N_kf\|_{L^2(\mu)}\, \|D^N_jg\|_{L^2(\mu)}.
\end{eqnarray*}
Since the matrix $\{2^{-\nu|j-k|}\}_{j,k}$ originates an operator bounded on
$\ell^2$, we obtain
\begin{eqnarray*}
|\langle T_{J,K}f,\, g\rangle| & \leq &
C\,\biggl(\sum_k \|D^N_kf\|_{L^2(\mu)}^2 \biggr)^{1/2} \,
\biggl(\sum_j\|D^N_jg\|_{L^2(\mu)}^2 \biggr)^{1/2}\\
& \leq & C\,N^2 \,\|f\|_{L^2(\mu)}\,\|g\|_{L^2(\mu)}.
\end{eqnarray*}
\end{proof}

\begin{lemma} \label{convernormal}
For $f,g \in \CC^\infty(\R^d)$ with compact support, we have
$$\lim_{m\to+\infty}\,
\Bigl\langle {\textstyle
T\Bigl(\sum_{|k|\leq m} D_k\, D_k^N f\Bigr),\,
\sum_{|j|\leq m} D_j\, D_j^N g}
\Bigr\rangle = \bigl\langle T(\Phi_N f),\,\Phi_N g\bigr\rangle.$$
\end{lemma}

\begin{proof} We know that $P_mf :=\sum_{|j|\leq m} D_j\, D_j^N f$ and
$P_mg := \sum_{|k|\leq m} D_k\, D_k^N g$ converge respectively to $\Phi_Nf$
and $\Phi_N g$ in $L^2(\mu)$.
Since we are assuming that the kernel $k(x,y)$ of $T$ is
bounded, we have $\|k(x,\cdot)\|_{L^p(\mu)},
\,\|k(\cdot,y)\|_{L^p(\mu)}\leq C$, for all $x,y$ and $1<p\leq\infty$.
As a consequence,
$T(P_m f)$ converges to $T(\Phi_N f)$ uniformly on $\R^d$ as $m\to\infty$.
Therefore, for any compact set $E\subset\R^d$, we have
$$\lim_{m\to\infty} \int_E T(P_m f)\,P_mg\,d\mu =
 \int_E T(\Phi_N f)\,\Phi_Ng\,d\mu.$$

It can be checked that there exists some constant $C_8$ independent of
$m$ such that the kernels $p_m(x,y)$ of the operators $P_m$ satisfy
the inequality
\begin{equation} \label{pm22}
|p_m(x,y)| \leq \frac{C_8}{|x-y|^n}.
\end{equation}
This an easy estimate that is left to reader.

We take $R>0$ so that $\supp(f),\, \supp(g) \subset B(0,R)$,
and $x_0$ with $|x_0|\geq10R$. By \rf{pm22}
we have
\begin{equation} \label{pmab}
|P_mf(y)|,\,|P_mg(y)|  \leq \frac{C}{|y|^n}\qquad \mbox{if $|y|\geq10R$},
\end{equation}
where $C$ may depend on $f$ and $g$ and, in particular, we may have $y=x_0$.
We split $|T(P_mf)(x_0)|$ as follows
\begin{eqnarray*}
|T(P_mf)(x_0)|\! & \leq & \!|T[(P_mf)\,\chi_{B(0,|x_0|/2)}] (x_0)|
+ |T[(P_mf)\,\chi_{\R^d\setminus B(0,|x_0|/2)}] (x_0)| \\
& = &\! A+ B.
\end{eqnarray*}
Let us estimate $A$:
\begin{eqnarray*}
A & \leq & \int_{ B(0,|x_0|/2)} C_{8}\,\frac{|P_mf(y)|}{|x_0-y|^n}\, d\mu(y)\\
& \leq & C\,\|P_mf\|_{L^2(\mu)} \left(\int_{ B(0,|x_0|/2)}
\frac{1}{|x_0-y|^{2n}}\, d\mu(y)\right)^{1/2} \\
& \leq & C\,\|P_mf\|_{L^2(\mu)}\frac{1}{|x_0|^{n/2}} \, = \,
\frac{C}{|x_0|^{n/2}}.
\end{eqnarray*}
Let us consider the term $B$. Since $k(x,\cdot)$ is in $L^2(\mu)$
uniformly on $x$, we have
$B \leq C\,\|(P_mf)\,\chi_{\R^d\setminus B(0,|x_0|/2)}\|_{L^2(\mu)}.$
From \rf{pmab} we get
$$\|(P_mf)\,\chi_{\R^d\setminus B(0,|x_0|/2)}\|_{L^2(\mu)} \leq
\frac{C}{|x_0|^{n/2}}.$$
Therefore, $|T(P_mf)(x_0)| \leq C/|x_0|^{n/2}$ (for $|x_0|\geq10R$).

Now we write
$$\int  T(P_m f)\,P_mg\,d\mu = \int_{B(0,10R)}\!\!
T(P_m f)\,P_mg\,d\mu +
\int_{\R^d\setminus B(0,10R)}\!\!  T(P_m f)\,P_mg\,d\mu.$$
The first integral on the right hand side tends to
$\int_{B(0,10R)} T(\Phi_N f)\,\Phi_N g\,d\mu$. The second one tends
to $\int_{\R^d\setminus B(0,10R)} T(\Phi_N f)\,\Phi_N g\,d\mu$, by an
application of the dominated convergence theorem, because
$|T(P_mf)(x_0)\,P_mg(x_0)|\leq C/|x_0|^{3n/2}$ if $x_0\in \R^d\setminus B(0,10R).$
\end{proof}

\vv
\noindent{\em Proof of Lemma \ref{t10}.} From the last lemmas we get
$$|\langle T\,\Phi_N f,\, \Phi_Ng\rangle | \leq C\, \|f\|_{L^2(\mu)}\,
\|g\|_{L^2(\mu)}.$$
That is, $\Phi_N^*\,T\,\Phi_N$ is bounded on $L^2(\mu)$, which implies that
$T$ is bounded on $L^2(\mu)$, since $\Phi_N^{-1}$ exists and is bounded.
\qed


\subsection{The proof of Lemma \ref{difi}}  \label{provadifi}

In next lemma we recall (without proof) a well known estimate that we
will need.

\begin{lemma}  \label{triv}
Let $\vphi,\,\psi$ be $L^1(\mu)$ functions supported on cubes $Q$ and $R$
respectively, with $\dist(Q,R)>\ell(Q)/2$.
If $\int \vphi \,d\mu=0$, then
\begin{equation}  \label{estimbasic}
|\langle T\vphi,\,\psi\rangle| \leq
C\,\frac{\ell(Q)^\delta}{\dist(Q,R)^{n+\delta}}\,\|\vphi\|_{L^1(\mu)}\,
\|\psi\|_{L^1(\mu)}.
\end{equation}
\end{lemma}

To prove Lemma \ref{difi}, we will change the
notation. We set $Q_i = Q_{x,i}$ and $R_i= Q_{y,i}$ for all $i$.
Also, we write $\vphi=u_{x,j}$ and $\psi=u_{y,k}$. Thus $\vphi$ is supported
on $Q_{j-2}$ and $\psi$ on $R_{k-2}$. We denote the centers of these cubes by
$x_0$ and $y_0$ respectively.
So in the case $2Q_{j-3}\cap2R_{k-3}=\varnothing$
(that is (a) in Lemma \ref{difi}), we have to prove that
\begin{equation}  \label{cvb0}
|\langle T\vphi,\, \psi\rangle| \leq
C\,2^{-\nu|j-k|}\,\frac{(\ell(Q_{j-2}) \wedge \ell(R_{k-2}))^{\delta/2}}{
(\ell(Q_{j-2}) + \ell(R_{k-2}) + |x_0-y_0|)^{n+\delta/2}},
\end{equation}
and if $2Q_{j-3}\cap2R_{k-3}\neq\varnothing$
(that is (b) in Lemma \ref{difi}), we have to show that
\begin{eqnarray} \label{cvb1}
|\langle T\vphi,\, \psi\rangle| & \leq &
\frac{C\,2^{-\nu|j-k|}}{(\ell(Q_j) + |x_0-y_0|)^n}\,
\chi_{Q_{j-7}}(y_0) \nonumber \\
&& \mbox{} +
\frac{C\,2^{-\nu|j-k|}}{(\ell(R_k) + |x_0-y_0|)^n}\, \chi_{R_{k-7}}(x_0).
\end{eqnarray}


\vspace{5mm}
\noi{\bf Proof of \rf{cvb0} for $2Q_{j-3}\cap2R_{k-3}=\varnothing$}.

Assume $j\geq k$, for example.
We have $Q_{k-2}\cap R_{k-2}=\varnothing$, because otherwise
$Q_{k-2}\subset R_{k-3}$, which implies $Q_{j-3}\cap R_{k-3}\neq\varnothing$.
Thus we get $|x_0-y_0|\geq \ell(Q_{k-2})/2$.

Now \rf{cvb1} is a direct consequence of Lemma \ref{triv}:
\begin{eqnarray*}
|\langle T\vphi,\,\psi\rangle| & \leq &
C\,\frac{(\ell(Q_{j-2}) \wedge \ell(R_{k-2}))^\delta}{(\ell(Q_{j-2}) +
\ell(R_{k-2})
+ |x_0-y_0|)^{n+\delta}} \,\|\vphi\|_{L^1(\mu)}\, \|\psi\|_{L^1(\mu)}\\
& \leq & C\,\frac{\ell(Q_{j-2})^{\delta/2}}{\ell(Q_{k-2})^{\delta/2}} \cdot
\frac{(\ell(Q_{j-2}) \wedge \ell(R_{k-2}))^{\delta/2}}{(\ell(Q_{j-2}) +
\ell(R_{k-2})
+ |x_0-y_0|)^{n+\delta/2}} \\
& \leq & C\, 2^{-\nu|j-k|}\,\frac{(\ell(Q_{j-2}) \wedge
  \ell(R_{k-2}))^{\delta/2}}{(\ell(Q_{j-2}) + \ell(R_{k-2})
+ |x_0-y_0|)^{n+\delta/2}}.
\end{eqnarray*}
Notice that we have used $|x_0-y_0|\geq \ell(Q_{k-2})/2$ in the second
inequality.
\qed


\vspace{5mm}
\noi{\bf Proof of \rf{cvb1} in the case $|j-k|>3,\,\,2Q_{j-3}\cap2R_{k-3}\neq
\varnothing$}.

We assume $j>k+3$.
Then we have $2Q_{j-3}\subset Q_{k-3} \subset R_{k-4}$.
If $x_0\in R_k$, then $\ell(Q_{j-2})\leq 2^{-\eta|(j-2)-k|}\,
\ell(R_k)\ll\ell(R_k)$, and so $Q_{j-2}\subset 2R_k$.

If $x_0\not\in R_k$, then $\ell(Q_{k+1})\leq 4|x_0-y_0|$ (otherwise
$R_k\subset Q_{k+1}$, which is not possible). Therefore
$\ell(Q_{j-2})\ll \ell(Q_{k+1}) \leq 4|x_0-y_0|$. So if we let $m\geq1$ be the
smallest integer such that $x_0\in 2^mR_k$, we will have
$$Q_{j-2} \subset 2^{m+1}R_k\setminus 2^{m-2}R_k$$
and
\begin{equation}  \label{difi2}
\dist(Q_{j-2},\,\R^d\setminus(2^{m+1}R_k\setminus 2^{m-2}R_k)) \approx
\ell(2^mR_k).
\end{equation}
Thus, in any case it is enough to prove  that
\begin{equation} \label{difi1}
|\langle T\vphi,\,\psi\rangle | \leq C\,\frac{2^{-\nu|j-k|}}{\ell(2^mR_k)^n}.
\end{equation}
(we take $m=0$ if $x_0\in R_k$).

We denote $L_m=2^{m+1}R_k\setminus 2^{m-2}R_k$. By Lemma \ref{triv}
(or a slight variant of it) and \rf{difi2}, we have
$$
|\langle T\vphi,\,(1-\chi_{L_m}) \psi\rangle|
\leq  C\,\frac{\ell(Q_{j-2})^\delta}{\dist(Q_{j-2},\,
\R^d\setminus L_m)^{n+\delta}}
\leq C\,\frac{\ell(Q_{j-2})^\delta}{\ell(2^mR_k)^{n+\delta}}.
$$
Arguing as above, we get
\begin{eqnarray*}
\ell(Q_{j-2}) & \leq & C\,2^{\eta|j-k|}\,\ell(Q_{k+1})\\
& \leq & C\,2^{-\eta|j-k|}\, (\ell(R_k) + |x_0-y_0|) \,
\leq \, C\,2^{-\eta|j-k|}\, \ell(2^mR_k).
\end{eqnarray*}
Therefore,
\begin{equation}\label{difi01}
\frac{\ell(Q_{j-2})^\delta}{\ell(2^mR_k)^{n+\delta}} \leq
C\, \frac{2^{-\nu|j-k|}}{\ell(2^mR_k)^n},
\end{equation}
and so
\begin{equation}\label{difi00}
|\langle T\vphi,\,(1-\chi_{L_m}) \psi\rangle| \leq
 C\,\frac{2^{-\nu|j-k|}}{\ell(2^mR_k)^n}.
\end{equation}

We have to prove that
\rf{difi1} holds for $\psi\,\chi_{L_m}$.
Since $T1\equiv0$, we have
\begin{equation}  \label{difi5}
\langle T\vphi,\,\psi\,\chi_{L_m}\rangle = \langle T\vphi,\,(\psi -
\psi(x_0))\,\chi_{L_m}\rangle -
\psi(x_0)\,\langle T\vphi, \chi_{\R^d\setminus L_m}\rangle = A + B.
\end{equation}
Taking into account $\int \vphi\,d\mu = 0$ and \rf{difi2}, by standard
estimates (similar to the ones of Lemma \ref{triv}) we get
$$|\langle T\vphi, \chi_{\R^d\setminus L_m}\rangle| \leq C\,
\frac{\ell(Q_{j-2})^\delta}{\ell(2^mR_k)^{\delta}}.$$
Since
\begin{equation} \label{difi4}
\|\psi\,\chi_{L_m}\|_{L^\infty(\mu)}\leq \frac{C}{\ell(2^mR_k)^n},
\end{equation}
arguing as in \rf{difi01}, we obtain
\begin{equation} \label{difib}
|B|\leq  C\,\frac{\ell(Q_{j-2})^\delta}{\ell(2^mR_k)^{n+\delta}}
\leq C\,\frac{2^{-\nu|j-k|}}{\ell(2^mR_k)^n}.
\end{equation}

Now we will estimate the term $A$ in \rf{difi5}.
We consider a bump function $w$ such that
$\chi_{2Q_{j-2}} \leq w \leq \chi_{4Q_{j-2}},$
with $|w'|\leq C/\ell(Q_{j-2})$.
We write
\begin{equation} \label{difia}
A =\bigl\langle T\vphi, (\psi -\psi(x_0))\,(\chi_{L_m}- w^2) \bigr\rangle
+ \bigl\langle T\vphi, (\psi -\psi(x_0))\, w^2 \bigr\rangle = A_1 + A_2.
\end{equation}
Since $\int \vphi\,d\mu =0$, we have
\begin{equation} \label{difi6}
|A_1| \leq \int_{y\in L_m\setminus 2Q_{j-2}} \int_{x\in Q_{j-2}}
\frac{C\,|x-x_0|^\delta}{|y-x_0|^{n+\delta}}\,|\vphi(x)|\,|\psi(y)-\psi(x_0)|
\,d\mu(x)\,d\mu(y).
\end{equation}
Recall that, by Lemma \ref{losire},
\begin{equation}  \label{difi7}
|\psi(y) - \psi(x)| \leq C\,\frac{|y-x|}{\ell(Q_k)}\cdot
\frac{1}{\ell(2^mR_k)^n} \leq C\,\frac{|y-x|^\delta}{\ell(Q_k)^\delta}\cdot
\frac{1}{\ell(2^mR_k)^n}
\end{equation}
if $y,x\in Q_k$.
We would like to plug this estimate (with $x=x_0$) into \rf{difi6}.
However we don't know if
\rf{difi7} holds for $y\in L_m\setminus Q_k$. Thus we split the double
integral in \rf{difi6} into two pieces:
$$|A_1| \leq \int_{y\in (L_m\cap Q_k)\setminus Q_{j-2}} \int_{x\in Q_{j-2}}
+ \int_{y\in L_m\setminus Q_k} \int_{x\in Q_{j-2}} = A_{1,1} + A_{1,2}.$$

We consider the integral $A_{1,1}$ first.
Using \rf{difi7}, we obtain:
$$A_{1,1} \leq \frac{C\,\ell(Q_{j-2})^\delta}{\ell(Q_k)^\delta\,
\ell(2^mR_k)^{n}}
\int_{x\in Q_{j-2}} |\vphi(x)| \int_{y\in Q_k\setminus 2Q_{j-2}}
\frac{1}{|y-x_0|^{n}}\,d\mu(y) \,d\mu(x).$$
Observe that
\begin{eqnarray*}
\int_{y\in Q_k\setminus 2Q_{j-2}}
\frac{1}{|y-x_0|^{n}}\,d\mu(y)
& \leq &
C\int_{y\in Q_k\setminus 2Q_{j-2}}
\frac{\ell(Q_k)^{\delta/2}}{|y-x_0|^{n+\delta/2}}\,d\mu(y) \\
&\leq & C\,\frac{\ell(Q_k)^{\delta/2}}{\ell(Q_{j-2})^{\delta/2}}.
\end{eqnarray*}
Therefore,
\begin{equation}  \label{difia11}
A_{1,1} \leq C\,\frac{\ell(Q_{j-2})^{\delta/2}}{\ell(Q_k)^{\delta/2}\,
\ell(2^mR_k)^{n}}
\,\|\vphi\|_{L^1(\mu)}
\leq C\, \frac{2^{-\nu|j-k|}}{\ell(2^mR_k)^n}.
\end{equation}

Let us consider $A_{1,2}$ now. Using \rf{difi4} we get
\begin{eqnarray} \label{difia12}
A_{1,2} & = & \int_{y\in L_m\setminus Q_k} \int_{x\in Q_{j-2}}
\frac{C\,|x-x_0|^\delta}{|y-x_0|^{n+\delta}}\,|\vphi(x)|\,|\psi(y)-\psi(x_0)|
\,d\mu(x)\,d\mu(y) \nonumber \\
& \leq & C\,\frac{\ell(Q_{j-2})^\delta}{\ell(2^mR_k)^n}
\int_{x\in Q_{j-2}} |\vphi(x)|
\int_{y\in L_m\setminus Q_k}\frac{1}{|y-x_0|^{n+\delta}}\,d\mu(y)\,d\mu(x)
\nonumber \\
& \leq &  C\,\frac{\ell(Q_{j-2})^\delta}{\ell(2^mR_k)^n \,\ell(Q_k)^\delta}
\int_{y\in L_m\setminus Q_k}\frac{1}{|y-x_0|^{n}}\,d\mu(y)\nonumber \\
& \leq &  C\,\frac{2^{-\nu|j-k|}}{\ell(2^mR_k)^n}\,(1+\delta(Q_k,R_{k-2}))
\, \leq \,  C\,\frac{2^{-\nu|j-k|}}{\ell(2^mR_k)^n}.
\end{eqnarray}

It only remains to estimate the term $A_2$ in \rf{difia}. As in \cite{DJS},
we introduce the term
$A_2' = \langle T(\vphi\,(\psi -\psi(x_0))\,w),\, w \rangle.$
First we will estimate the difference $|A_2-A_2'|$, and later $A_2'$.

We write $\psi_0 = (\psi -\psi(x_0))\,w$, and then we have
$$A_2 - A_2' = \iint (\psi_0(y)-\psi_0(x))\, k(x,y)\,
\vphi(x)\,w(y)\,d\mu(x)\,d\mu(y).$$
Thus
$$|A_2 - A_2'| \leq C\,\|\psi_0\|_{lip1} \iint  \frac{1}{|x-y|^{n-1}}
\, |\vphi(x)|\,w(y)\,d\mu(x)\,d\mu(y).$$
Since $\int_{4Q_{j-2}} |x-y|^{1-n}\,d\mu(x) \leq C\,\ell(Q_{j-2})$
for any $y\in 4Q_{j-2}$, we obtain
$$|A_2 - A_2'| \leq C\,\|\psi_0\|_{lip1} \,\ell(Q_{j-2}) \int |\vphi(x)|
\,d\mu(x) \leq C\,\|\psi_0\|_{lip1} \,\ell(Q_{j-2}).$$
For $x,y\in 4Q_{j-2}\subset Q_{k}$ \rf{difi7} holds, and so
\begin{eqnarray*}
\|\psi_0\|_{lip1} & \leq & \|w\|_\infty\, \|\psi - \psi(x_0)\|_{lip1,4Q_{j-2}}
+ \|w\|_{lip1}\,\|\psi - \psi(x_0)\|_\infty \\
& \leq & \frac{C}{\ell(Q_k)\,\ell(2^mR_k)^n}.
\end{eqnarray*}
So we get
\begin{equation}  \label{difidif}
|A_2 - A_2'| \leq C\,\frac{\ell(Q_{j-2})}{\ell(Q_k)\,\ell(2^mR_k)^n}\leq
C\,\frac{2^{-\nu|j-k|}}{\ell(2^mR_k)^n}.
\end{equation}

Finally we have to deal with $A_2'$. We write
$$A_2' = \langle \vphi\,\psi_0,\, T^*(\chi_{2Q_{j-2}}) \rangle
+ \langle \vphi\,\psi_0,\, T^*(w-\chi_{2Q_{j-2}}) \rangle.$$
Since $T^*$ is bounded on $L^p(\mu)$ over characteristic
functions of cubes for $p=2$, and for $p=n/(n-1)$ in the case
$n>1$, and it is also bounded from
$L^p(\mu\mid \R^d\setminus 2Q_{j-2})$ into $L^p(\mu\mid Q_{j-2})$
(for any $p\in (1,\infty)$), we obtain
$|A_2'| \leq C\, \|\vphi\,\psi_0\|_{L^p(\mu)}\, \mu(4Q_{j-2})^{1/p'}.$
Now we can estimate the $L^p(\mu)$ norm of $\vphi\,\psi_0$ using \rf{difi7}:
\begin{eqnarray}  \label{wertq}
\|\vphi\,\psi_0\|_{L^p(\mu)}^p & \leq &
\int_{Q_{j-2}}\left[\frac{|x-x_0|}{(\ell(Q_j) + |x-x_0|)^n\, \ell(Q_k)\,
\ell(2^mR_k)^n}\right]^p\, d\mu(x) \nonumber \\
& \leq & \frac{C}{\left[\ell(Q_k)\, \ell(2^mR_k)^n\right]^p}
\left( \int_{Q_j} \frac{1}{\ell(Q_j)^{(n-1)p}}\, d\mu(x) \right. \nonumber \\
&& \mbox{} +\left.
\int_{Q_{j-2}\setminus Q_j} \frac{1}{|x-x_0|^{(n-1)p}}\, d\mu(x)\right).
\end{eqnarray}
If $n>1$, we choose $p=n/(n-1)$, and we get
$$\|\vphi\,\psi_0\|_{L^p(\mu)} \leq \frac{C}{\ell(Q_k)\, \ell(2^mR_k)^n}.$$
Therefore, we have
\begin{equation}  \label{difi20}
|A_2'| \leq C\,\frac{\ell(Q_{j-2})^{n/p'}}{\ell(Q_k)\, \ell(2^mR_k)^n}
= C\,\frac{\ell(Q_{j-2})}{\ell(Q_k)\, \ell(2^mR_k)^n}
\leq C\,\frac{2^{-\nu|j-k|}}{\ell(2^mR_k)^n}.
\end{equation}
If $n\leq1$, we take $p=2$, and from \rf{wertq} we obtain
$$\|\vphi\,\psi_0\|_{L^2(\mu)} \leq
C\,\frac{\ell(Q_{j-2})^{1-n/2}}{\ell(Q_k)\, \ell(2^mR_k)^n}.$$
Then we also have
\begin{equation}  \label{difi20'}
|A_2'| \leq  C\,\frac{\ell(Q_{j-2})}{\ell(Q_k)\, \ell(2^mR_k)^n}
\leq C\,\frac{2^{-\nu|j-k|}}{\ell(2^mR_k)^n}.
\end{equation}

From \rf{difi00}, \rf{difib}, \rf{difia11}, \rf{difia12}, \rf{difidif},
\rf{difi20}, and \rf{difi20'} we obtain \rf{difi1}.
\qed


\vspace{5mm}
\noi{\bf Proof of \rf{cvb1} in the case
  $|j-k|\leq3,\,\,2Q_{j-3}\cap2R_{k-3}
\neq \varnothing$}.

Observe that in this case, since $|j-k|\leq3$,
then $Q_{j-2} \subset R_{k-7}$ and $R_{k-3}\subset Q_{j-7}$.
Assume first that $10Q_j\cap 10R_k=\varnothing$. We denote $d=|x_0-y_0|$ and
$A= Q_{j-7}\setminus B(y_0,d/4) $. Notice that $d\approx
\dist(Q_j,R_k)$. We will show that
\begin{equation}  \label{difi30}
|\langle T\vphi,\,\psi\rangle| \leq \frac{C}{d^n},
\end{equation}
which implies \rf{cvb1}.

Let $w_A$ be a $\CC^1$ function such that $0\leq w_A\leq1$, $w_A\equiv1$ on $A$,
$w_A\equiv0$ on $\R^d\setminus U_{d/20}(A)$ (where $U_\ve(A)$ is the
$\ve$-neighbourhood of $A$), with $|w_A'|\leq C/d$. We split
$\langle T\vphi,\,\psi\rangle$ as follows:
\begin{equation}  \label{difi30.5}
\langle T\vphi,\,\psi\rangle =
\langle T\vphi,\,\psi w_A\rangle  +  \langle T\vphi,\,\psi(1-w_A)\rangle
= I + J.
\end{equation}

First we will estimate $I$. Observe that
\begin{equation}  \label{difi31}
|\psi\,w_A| \leq \frac{C}{d^n}
\end{equation}
(this inequality will be basic in our arguments). Let us consider the term
$I'=\langle T(\vphi\,\psi\,w_A),\,\chi_{9Q_{j-7}}\rangle$. We have
\begin{eqnarray*}
\lefteqn{|I-I'| \, = \,|\langle T\vphi,\,\psi w_A\chi_{9Q_{j-7}}\rangle
-\langle T(\vphi\,\psi\,w_A),\,\chi_{9Q_{j-7}}\rangle}| && \\
& \leq & \iint \frac{C}{|x-y|^n}\, |\vphi(x)|\,|\psi(x)w_A(x) -
\psi(y)w_A(y)|\,\chi_{9Q_{j-7}}(y)\,d\mu(x)\,d\mu(y)\\
& = & \int_x \int_{y\in 9Q_{j-7}\cap Q_{x,k+1}}
+ \int_x \int_{y\in 9Q_{j-7}\setminus Q_{x,k+1}} =\, H_1 + H_2.
\end{eqnarray*}

Let us consider the integral $H_1$:
\begin{eqnarray*}
H_1 & \leq & C \int |\vphi(x)|\, \|\psi w_A\|_{lip_1,Q_{x,k+1}}
\int_{y\in Q_{x,k+1}} \frac{1}{|x-y|^{n-1}}\,d\mu(y)\, d\mu(x) \\
& \leq & C \int|\vphi(x)|\,\|\psi w_A\|_{lip_1,Q_{x,k+1}}\,\ell(Q_{x,k+1})\,
d\mu(x).
\end{eqnarray*}
By Lemma \ref{losire}, we have
$\|\psi \|_{lip_1,Q_{x,k+1}} \leq {C}/(\ell(Q_{x,k})\,d^n),$
and so
$$
\|\psi w_A\|_{lip_1,Q_{x,k+1}}  \,\leq \, \|\psi \|_{lip_1,Q_{x,k+1}} +
\|\psi \|_\infty \,\|w_A\|_{lip_1}
\,\leq \,\frac{C}{\ell(Q_{x,k})\,d^n} + \frac{C}{d^{n+1}}.
$$
Since $\ell(Q_{x,k+1}) \leq 10d$ (otherwise $Q_{x,k+1}\supset R_k$, which is not
possible), we obtain
$$H_1  \leq  \frac{C}{d^n} \int |\vphi(x)|\,d\mu(x)\leq  \frac{C}{d^n}.$$

Let us turn our attention to $H_2$. From \rf{difi31} we get
\begin{eqnarray*}
H_2 & \leq & \frac{C}{d^n} \int |\vphi(x)| \int_{y\in 9Q_{j-7}\setminus
Q_{x,k+1}}  \frac{1}{|x-y|^{n}}\, d\mu(y) \,d\mu(x) \\
& \leq & \frac{C}{d^n}\int |\vphi(x)| \,\delta(Q_{x,k+1},9Q_{j-7})\,d\mu(x)
\, \leq \, \frac{C}{d^n}.
\end{eqnarray*}

Now we will estimate $I'$. We consider the annuli $C_i= 3^iQ_j\setminus
3^{i-1}Q_j$ ($i\geq1$), $C_0=Q_j$, and some neighbourhoods of them $\wt{C}_i
=3^{i+1}Q_j\setminus 3^{i-2}Q_j$ ($i\geq1$), $\wt{C}_0=3Q_j$. We write
\begin{eqnarray*}
I' & = & \sum_{i=0}^{N_0} \langle T(\vphi\psi w_A\chi_{C_i}),\,
\chi_{9Q_{j-7}}\rangle \\
& = & \sum_{i=0}^{N_0} \langle T(\vphi\psi w_A\chi_{C_i}),\,
\chi_{9Q_{j-7}\cap\wt{C}_i} \rangle +
\sum_{i=0}^{N_0} \langle T(\vphi\psi w_A\chi_{C_i}),\,
\chi_{9Q_{j-7}\setminus\wt{C}_i} \rangle \\
& = & \sum_{i=0}^{N_0} I'_{1,i} + \sum_{i=0}^{N_0} I'_{2,i}.
\end{eqnarray*}
From the $L^2(\mu)$ boundedness of $T^*$ over characteristic functions of
cubes (as shown in Lemma \ref{caracterist}) we get
$$
|I'_{1,i}| \leq C\,\|\vphi\psi w_A\chi_{C_i}\|_{L^2(\mu)}\,
\mu(2\wt{C}_i)^{1/2} \leq  \frac{C}{d^n} \, \|\vphi\chi_{C_i}\|_{L^2(\mu)}\,
\mu(2\wt{C}_i)^{1/2}.
$$
Since $\|\vphi\chi_{C_i}\|_{L^2(\mu)} \leq C\mu(C_i)^{1/2}/3^i\ell(Q_j)^n$, we
obtain
$$\sum_{i=0}^{N_0}|I'_{1,i}| \leq \frac{C}{d^n}
\sum_{i=0}^{N_0}\frac{\mu(2\wt{C}_i)}{3^i\ell(Q_j)^n} \leq
\frac{C}{d^n}\,(C + \delta(Q_j,9Q_{j-7})) \leq \frac{C}{d^n}.$$
Let us consider the terms $I'_{2,i}$. For
$y\not\in \wt{C}_i$ we have
\begin{eqnarray*}
|T(\vphi\psi w_A\chi_{C_i})(y)| &\leq & \frac{C}{(|y-x_0|+\ell(Q_j))^n}
 \, \|\vphi\psi
w_A\chi_{C_i})\|_{L^1(\mu)} \\
& \leq & \frac{C}{(|y-x_0|+\ell(Q_j))^n\,d^n}\, \|\vphi\chi_{C_i}\|_{L^1(\mu)}.
\end{eqnarray*}
Therefore
\begin{eqnarray*}
|I'_{2,i}| & \leq & \frac{C\,\|\vphi\chi_{C_i}\|_{L^1(\mu)}}{d^n}
\int_{9Q_{j-7}} \frac{1}{(|y-x_0|+\ell(Q_j))^n}\, d\mu(y) \\
& \leq & \frac{C\,\|\vphi\chi_{C_i}\|_{L^1(\mu)}}{d^n},
\end{eqnarray*}
and so
$\ds \sum_{i=0}^{N_0}|I'_{2,i}| \leq \frac{C}{d^n}.$
Thus we have shown that $|I|\leq C/d^n$.

Now we have to consider the term
$J$ of \rf{difi30.5}. We take $B=R_{k-7}\setminus B(x_0,d/4)$, and we let $w_B$
be a $\CC^1$ function such that $0\leq w_B\leq1$, $w_B\equiv1$ on $B$,
$w_B\equiv0$ on $\R^d\setminus U_{d/20}(B)$, and $|w_B'|\leq C/d$. We write
$$J = \langle T(\vphi w_B),\, \psi\,(1-w_A)\rangle +
\langle T(\vphi \,(1-w_B)),\, \psi\,(1-w_A)\rangle = J_1 + J_2.$$
Now we have $|\vphi \,w_B| \leq C/d^n$. So the estimates for the
term $J_1$ are analogous to the ones for the term $I$ of \rf{difi30.5}.
We only have to interchange the roles of $\psi w_A$ and $\vphi w_B$,
$T$ and $T^*$, etc., and then we will get
$|J_1|\leq C/d^n$ too. The details are left to the reader.

Finally, we only have to deal with the term $J_2$. The estimates for this case
are straightforward. Since
$$\dist(\supp(\psi(1-w_A)),\, \supp(\vphi(1-w_B))) \geq C^{-1}\,d,$$
for $x\in \supp(\psi(1-w_A))$ we obtain
$$|T(\vphi(1-w_B))(x)| \leq \frac{C}{d^n}\, \|\vphi(1-w_B))\|_{L^1(\mu)}
\leq \frac{C}{d^n}.$$
Thus $|J_2| \leq C\, \|\psi(1-w_A))\|_{L^1(\mu)}/d^n
\leq C/d^n.$
Therefore, \rf{difi30} holds if $10Q_j\cap 10R_k=\varnothing$.

\vv
The case $10Q_j\cap 10R_k\neq\varnothing$ is simpler. Assume for example
$\ell(Q_j) \leq \ell(R_k)$. Then it is enough to show that
\begin{equation} \label{difi40}
|\langle T\vphi,\, \psi\rangle| \leq \frac{C}{\ell(R_k)^n}.
\end{equation}
Notice that we have
$\|\psi\|_{L^\infty(\mu)} \leq C/\ell(R_k)^n$. Then in this case it is not
necessary to split $\langle T\vphi,\, \psi\rangle$ into two terms $I$ and $J$
as in \rf{difi30}. Estimates similar to the ones used for the term $I$
will yield \rf{difi40}. We omit the detailed arguments again.
\qed


\section{The paraproduct}  \label{sec8}

For technical reasons, we need to introduce a class of operators slightly
larger than the class of CZO's that we have considered.

\begin{definition} \label{defhczo}
We say that $k(x,y)$ is a H\"ormander-Calder\'on-Zygmund kernel if
\begin{itemize}
\item[(1)] $\ds |k(x,y)| \leq \frac{C_1}{|x-y|^n}$ if $x\neq y$,
\item[(2')] for any $x,x'\in\supp(\mu)$,
$$ \int_{|y-x|\geq2|x-x'|}\bigl(|k(x,y)-k(x',y)| + |k(y,x)-k(y,x')|\bigr)
\,d\mu(y) \leq C_2'.$$
\end{itemize}
We say that $T$ is a H\"ormander-Calder\'on-Zygmund operator (HCZO)
if $T$ is associated to the kernel $k(x,y)$ as in \rf{eq***}.
\end{definition}

Condition (2') in the definition above is called H\"ormander's condition.

Recall that in Section \ref{lipa} we have defined
$\Phi_N = \sum_{k\in\Z} D_j^N \,D_j$ and in Lemma \ref{lemadk} we have shown
that $\Phi_N\to I$ as $N\to\infty$ in the operator norm of $L^2(\mu)$, and so
$\Phi_N$ is invertible in $L^2(\mu)$ for $N$ big enough. We will show
analogous results for $\Phi_N$ on $L^p(\mu)$, $1<p<\infty$, and $\rbmo(\mu)$.
We need the following lemma.

\begin{lemma} \label{rbmorbmo}
If $T$ is a HCZO bounded on $L^2(\mu)$ with $T(1)=0$, then $T$ can be extended
boundedly from $\rbmo(\mu)$ into $\rbmo(\mu)$.
\end{lemma}

In \cite[Theorem 2.11]{Tolsa3} it is shown that if $T$ is a CZO which
is bounded on $L^2(\mu)$, then it is also bounded from $L^\infty(\mu)$ into
$\rbmo(\mu)$. A small modification of these arguments yields the result
stated in Lemma \ref{rbmorbmo}, as in the usual
doubling case (see \cite[Lemma 2.7]{DJS}, for example).

\begin{definition} \label{normhczo}
Let $T$ be a HCZO bounded on $L^2(\mu)$. We define
the HCZO norm of $T$ as $\|T\|_{HCZO} :=
\|T\|_{2,2} + C_1 + C_2'$, where $C_1$ and $C_2'$
are the best constants that appear in the conditions (1) and (2') defining a
HCZ kernel. Moreover, we say that a sequence of
linear operators $\{T_k\}_k$ converges to some linear operator $T$ in HCZO
norm if $\|T-T_k\|_{HCZO}\to 0$ as $k\to\infty$.
\end{definition}

\begin{lemma}  \label{inverbmo}
The operator $\Phi_N$ tends to $I$ in HCZO norm as $N\to\infty$.
Moreover, $\Phi_N$ can be extended boundedly
on $L^p(\mu)$, $1<p<\infty$, and from $\rbmo(\mu)$ into
$\rbmo(\mu)$. For $N$ big enough it is invertible in $L^p(\mu)$ (with $N$
depending on $p$) and in $\rbmo(\mu)$.
\end{lemma}

\begin{proof}
Notice that we only have to show that $I-\Phi_N$ is a HCZO such that
$\|I-\Phi_N\|_{HCZO}\to0$ as $N\to\infty$, taking into account Lemma
\ref{rbmorbmo} and the fact that HCZO's are bounded
on $L^p(\mu)$ (see \cite{NTV2}, and also \cite{Tolsa4} for a different proof).

We have already seen in Lemma \ref{lemadk} that $\Phi_N\to I$ as $N\to\infty$
in the operator norm of $L^2(\mu)$. Thus it only remains to
see that $I-\Phi_N$ is a HCZO and that the constants $C_1$ and $C_2'$ in
Definition \ref{defhczo} tend to $0$ as $N\to\infty$.
Recall that $I-\Phi_N = \sum_{|k|>N} E_k$, with $E_k = \sum_{j\in\Z} D_{j+k}
\,D_{j}$.

First we deal with the inequality (1) in Definition \ref{defhczo}.
In \rf{eqxxx} we have shown that if $k\geq2$, then the kernel $K_{j+k,j}(x,y)$
of $D_{j+k}\,D_j$ satisfies
\begin{equation} \label{eqxxx'}
|K_{j+k,j}(x,y)| \leq C\,2^{-\eta k}\,
\frac{1}{(\ell(Q_{x,j}) + \ell(Q_{y,j}) + |x-y|)^n}.
\end{equation}
Moreover, just below \rf{eqxxx} we have seen that
$K_{j+k,j}(x,y)=0$ if $y\not \in  Q_{x,j-3}$ or  $x\not \in  Q_{y,j-3}$.
For $x,y\in\supp(\mu)$, $x\neq y$,
let $j_0$ be the largest integer such that $y\in Q_{x,j_0}$. Since
$y\not\in Q_{x,j_0+h}$ for $h\geq1$, we get $K_{j+k,j}(x,y)=0$ if $j\geq
j_0+4$. Taking into account that for $j\leq j_0$ we have $|x-y|\leq
\ell(Q_{x,j_0}) \leq 2^{-\eta|j-j_0|}\,\ell(Q_{x,j})$, from \rf{eqxxx'}
it easily follows that
$$
\sum_{j\in\Z} |K_{j+k,j}(x,y)| \leq C\,2^{-\eta |k|}\,\frac{1}{|x-y|^n}.
$$
An analogous estimate can be obtained for $k\leq-2$.
Thus the kernel $K_N(x,y)$ of $I-\Phi_N$ satisfies
\begin{equation}  \label{eqwww1}
|K_N(x,y)|\leq C\,2^{-\eta N}\,\frac{1}{|x-y|^n}.
\end{equation}

Now we have to show that the constant $C_2'$ in Definition \ref{defhczo}
corresponding to the kernel of $I-\Phi_N$ tends to $0$ as $N\to\infty$.
First we will deal with the term $I:=|K_{j+k,j}(x,y)-K_{j+k,j}(x',y)|$,
assuming $k\geq N\geq 10$. Let $h_0$ be the largest integer such that $x'\in
Q_{x,h_0}$. Using \rf{eqwww1} it is easy to check that
$$\sum_{k\geq N} \sum_{j\in\Z} \int_{y\in Q_{x,h_0-10}\setminus B(x,2|x-x'|)}
I\,d\mu(y) \leq C\, 2^{-\eta N}.$$
So we only have to estimate the integral $\int_{\R^d\setminus Q_{x,h_0-10}}
I\,d\mu(y)$. Notice that
$\supp\, K_{j+k,j}(x,\cdot) \subset Q_{x,j-3}$ and
$\supp\, K_{j+k,j}(x',\cdot) \subset Q_{x',j-3}\subset Q_{x,j-4}
\cup Q_{x,h_0-10}$, and so
$\supp(I)\subset Q_{x,j-4}\cup Q_{x,h_0-10}$. Thus we may assume
$j-4\leq h_0-10$.
Let us consider the case $j+k>h_0$, that is, $x'\not\in Q_{x,j+k}$. By
\rf{eqwww21}, \rf{eqwww22} we obtain
$$\sum_{k\geq N}\,\, \sum_{
\substack{j:\,j-4\leq h_0-10\\ j>h_0-k}}
 \int I\,d\mu(y)
\leq C\,\sum_{k\geq N} 2^{-\eta k}\,(14+k) \leq C\,2^{-\eta N/2}.$$
Assume now $j+k\leq h_0$, that is $x'\in Q_{x,j+k}$ (and $k\geq N\geq 10$ too).
Observe that
\begin{multline}  \label{hc1}
K_{j+k,j}(x,y)-K_{j+k,j}(x',y) = \int (d_{j+k}(x,z) - d_{j+k}(x',z))\,
d_j(z,y) \,d\mu(z) \\
 =  \int (d_{j+k}(x,z) - d_{j+k}(x',z))\,
(d_j(z,y)-d_j(x,y)) \,d\mu(z).
\end{multline}
It is easily checked that the integrand above is null unless
$z\in Q_{x,j+k-3}$. Since $x'\in Q_{x,j+k}$, we have
$$|d_{j+k}(x,z) - d_{j+k}(x',z)| \leq C\,\frac{|x-x'|}{\ell(Q_{x,j+k})}\cdot
\frac{1}{(\ell(Q_{x,j+k}) + |x-z|)^n}.$$
Also, for $z\in Q_{x,j+k-3}\subset Q_{x,j}$,
\begin{equation} \label{hc2}
|d_j(z,y)-d_j(x,y)| \leq C\, \frac{|x-z|}{\ell(Q_{x,j})}\cdot
\frac{1}{(\ell(Q_{x,j}) + |x-y|)^n}.
\end{equation}
Therefore,
\begin{equation*}
\begin{split}
I & \leq C\, \frac{|x-x'|\,\ell(Q_{x,j+k-3})}{\ell(Q_{x,j})\, \ell(Q_{x,j+k})
\,(\ell(Q_{x,j}) + |x-y|)^n} \\
& \quad \times \int_{z\in Q_{x,j+k-3}}
\frac{1}{(\ell(Q_{x,j+k}) + |x-z|)^n}\,d\mu(z) \\
& \leq  C\, \frac{|x-x'|\,\ell(Q_{x,j+k-3})}{\ell(Q_{x,j})\, \ell(Q_{x,j+k})
\,(\ell(Q_{x,j}) + |x-y|)^n}.
\end{split}
\end{equation*}
Using $\ell(Q_{x,j+k-3})/\ell(Q_{x,j}) \leq C\,2^{-\eta k}$, we obtain
\begin{eqnarray*}
\int I \,d\mu(y)
& \leq & C\,2^{-\eta k} \frac{|x-x'|}{\ell(Q_{x,j+k})}
\int_{Q_{x,j-4}} \frac{1}{(\ell(Q_{x,j}) + |x-y|)^n}\,d\mu(y)\\
& \leq & C\,2^{-\eta k} \frac{|x-x'|}{\ell(Q_{x,j+k})}.
\end{eqnarray*}
Thus we get
\begin{multline*}
\sum_{k\geq N} \sum_{j:\,j+k\leq h_0}
\int |K_{j+k,j}(x,y)-K_{j+k,j}(x',y)|\,d\mu(y)\\  \leq
C \sum_{k\geq N} 2^{-\eta k} \sum_{j:\,j+k\leq h_0}
\frac{\ell(Q_{x,h_0})}{\ell(Q_{x,j+k})}
\leq C\sum_{k\geq N} 2^{-\eta k} \leq C\, 2^{-\eta N}.
\end{multline*}

Let us consider the term $J:=|K_{j+k,j}(y,x)-K_{j+k,j}(y,x')|$ now. As in the
case of the term $I$, we have
$$\sum_{k\geq N} \sum_{j\in\Z} \int_{y\in Q_{x,h_0-10}\setminus B(x,2|x-x'|)}
J\,d\mu(y) \leq C\, 2^{-\eta N},$$
and we only have to consider the integral $\int_{\R^d \setminus Q_{x,h_0-10}}
J\,d\mu(y)$. Moreover, it is easily seen that we also have
$\supp(J)\subset Q_{x,j-4}\cup Q_{x,h_0-10}$ in this case.
Thus we may assume $j-4\leq h_0$
again. Operating as above, by \rf{eqwww21}, \rf{eqwww22} we obtain
$$\sum_{k\geq N} \,\,\sum_{\substack{
j:\,j-4\leq h_0-10\\ j>h_0-k}} \int J\,d\mu(y)
\leq C\,2^{-\eta N/2}.$$
Suppose now that $j+k\leq h_0$, that is, $x'\in Q_{x,j+k}$. We have
$$J \leq \int |d_{j+k}(y,z)\, (d_j(z,x) - d_j(z,x'))|\, d\mu(z).$$
Since $x'\in Q_{x,h_0} \subset Q_{x,j}$, we have
$$|d_j(z,x) - d_j(z,x')| \leq C\, \frac{|x-x'|}{\ell(Q_{x,j})} \cdot
\frac{1}{(\ell(Q_{x,j}) + |x-z|)^n}.$$
Thus
\begin{eqnarray*}
J & \leq &  C\, \frac{|x-x'|}{\ell(Q_{x,j})}
\int_{z\in Q_{x,j-4}} |d_{j+k}(y,z)|\,\frac{1}{(\ell(Q_{x,j})+|x-z|)^n}
\,d\mu(z) \\
& \leq & C\, \frac{|x-x'|}{\ell(Q_{x,j})}
\left(\int_{\begin{subarray}{l} z\in Q_{x,j-4}\\|y-z|\leq |x-y|/2
\end{subarray}}
+ \int_{\begin{subarray}{l} z\in Q_{x,j-4}\\|y-z|> |x-y|/2\end{subarray}}
\right)
= J_1 + J_2
\end{eqnarray*}
Let us estimate $J_1$:
\begin{eqnarray*}
J_1 & \leq &  C\, \frac{|x-x'|}{\ell(Q_{x,j})}
\int |d_{j+k}(y,z)|\,\frac{1}{(\ell(Q_{x,j}) + |x-y|)^n} \,d\mu(z)\\
& \leq &  C\, \frac{|x-x'|}{\ell(Q_{x,j})} \cdot
\frac{1}{(\ell(Q_{x,j}) + |x-y|)^n}.
\end{eqnarray*}
We consider $J_2$ now. On the one hand we have
\begin{eqnarray*}
J_2 & \leq & C\, \frac{|x-x'|}{\ell(Q_{x,j})}
\int_{z\in Q_{x,j-4}} \frac{1}{|x-y|^n}\cdot\frac{1}{(\ell(Q_{x,j})+|x-z|)^n}
\,d\mu(z)\\
& \leq &  C\, \frac{|x-x'|}{\ell(Q_{x,j})\,|x-y|^n}.
\end{eqnarray*}
On the other hand,
$$J_2 \leq C\, \frac{|x-x'|}{\ell(Q_{x,j})}
\int |d_{j+k}(y,z)|\,\frac{1}{\ell(Q_{x,j})^n}\,d\mu(z)
\leq  C\, \frac{|x-x'|}{\ell(Q_{x,j})^{n+1}}.$$
Thus we have
$$J_2 \leq  C\, \frac{|x-x'|}{\ell(Q_{x,j})} \cdot
\frac{1}{(\ell(Q_{x,j}) + |x-y|)^n}$$
in any case. Therefore,
$$\int J\,d\mu(y) \leq C\, \frac{|x-x'|}{\ell(Q_{x,j})}
\int_{y\in Q_{x,j-4}} \frac{1}{(\ell(Q_{x,j}) + |x-y|)^n}\, d\mu(y) \leq
 C\, \frac{\ell(Q_{x,h_0})}{\ell(Q_{x,j})},$$
and so
\begin{eqnarray*}
\lefteqn{
\sum_{k\geq N} \,\sum_{j:\,j+k\leq h_0}
\int |K_{j+k,j}(y,x)-K_{j+k,j}(y,x')|\,d\mu(y)}&& \\
& \leq &
C \sum_{k\geq N} \,\sum_{j:\,j+k\leq h_0}
\frac{\ell(Q_{x,h_0})}{\ell(Q_{x,j})} \\
& \leq & C \sum_{k\geq N} \,\sum_{j:\,j+k\leq h_0}
\frac{\ell(Q_{x,j+k})^{1/2}}{\ell(Q_{x,j})^{1/2}} \cdot
\frac{\ell(Q_{x,h_0})^{1/2}}{\ell(Q_{x,j})^{1/2}} \\
& \leq & C\sum_{k\geq N} 2^{-\eta k/2} \leq C\, 2^{-\eta N/2}.
\end{eqnarray*}

\vv
When $k$ is  negative ($k\leq-N$), we have analogous estimates.
As a consequence, the kernel of $I-\Phi_N$ satisfies H\"ormander's condition
with constant $C_2'\leq C\,2^{-\eta N/2}$, and we are done.
\end{proof}

In order to prove the $T(1)$ theorem in the general case, we will introduce a
paraproduct. Given a fixed function $b\in \rbmo(\mu)$, we denote by
$P_{k,b}^N$ the operator of pointwise multiplication by
$D_k^N\,\Phi_N^{-1}(b)$. Then, for each positive integer $m$, we set
$$U_{m,b} = \sum_{k=-m}^m D_k\, P_{k,b}^N\,S_k.$$
We will show that the operators $U_{m,b}$ are uniformly bounded on
$L^2(\mu)$. A weak limit $U_b$  of the sequence $\{U_{m,b}\}_m$ will
be our required paraproduct, which will satisfy $U_b (1) = 0$ and
$U_b^*(1) = 0$.

Before dealing with the $L^2(\mu)$ boundedness of the operators $U_{m,b}$, we
need a suitable discrete version of Carleson's imbedding theorem.
Given $E\subset\R^d$, we denote
$$\wt{E} = \{(x,k)\in E\times\Z:\, (Q_{x,k})^\circ\subset E\}.$$
Our discrete version of Carleson's result is the following:

\begin{lemma}
For each $k\in\Z$ let $a_k(\cdot)$ be some non negative function and
let $\nu_k$ be the measure given by $d\nu_k =
a_k\,d\mu$. We denote $\nu=\sum_{k\in\Z} \bar{\nu}_k$, where $\bar{\nu}_k$
stands for the measure $\nu_k$ `transported' to $\R^d\times \{k\}$ (that is,
for $A\subset \R^d \times\Z$, $\bar{\nu}_k(A) = \nu_k\{x:\,(x,k)\in A\}$).
If
\begin{equation} \label{assum}
\sum_{j=k-2}^\infty \nu_j(Q) \leq C_9\,\mu(Q)
\end{equation}
for any doubling cube $Q$ of the $k$th generation, then we have
\begin{itemize}
\item[(a)] If $E\subset\R^d$ is open, then $\nu(\wt{E}) \leq C\,C_9\,\mu(E)$.
\item[(b)] For all $f\in L^2(\mu)$,
$$\sum_{k\in\Z} \|S_kf\|_{L^2(\nu_k)}^2 \leq C\,C_9\, \|f\|_{L^2(\mu)}^2.$$
\end{itemize}
\end{lemma}

\begin{proof} Let us see that (a) holds. We may assume that
$E$ is bounded. For each $x\in E$ we
choose the biggest cube $Q_{x,k}\subset E$ (i.e. with $k$ minimal). By
Besicovitch's covering theorem, there exists a family cubes $Q_{x_i,k_i}$ with
finite overlap such that $E = \bigcup_i Q_{x_i,k_i}$.
Therefore we have
\begin{equation} \label{tyr1}
\wt{E} = \bigcup_i (Q_{x_i,k_i} \times \Z) \cap \wt{E}.
\end{equation}

Observe also that if $(x,k) \in Q_{x_i,k_i} \cap \wt{E}$, then $k>k_i-3$.
Otherwise $(Q_{x,k_i-3})^\circ \subset E$, and since $x\in Q_{x_i,k_i}\subset
Q_{x_i,k_i-1}$, then
$$Q_{x_i,k_i-1} \subset Q_{x,k_i-2}  \subset(Q_{x,k_i-3})^\circ
\subset E,$$
which contradicts the maximality of $Q_{x_i,k_i}$.
Therefore we get
\begin{eqnarray*}
\nu((Q_{x_i,k_i}\times \Z)\cap \wt{E}) & = &
\nu((Q_{x_i,k_i} \times [k_i-2,+\infty)) \cap \wt{E}) \\
& \leq & \sum_{j=k_i-2}^{+\infty} \nu_j(Q_{x_i,k_i}) \, \leq \, C_9\,
\mu(Q_{x_i,k_i}).
\end{eqnarray*}
By \rf{tyr1} and the finite overlap of the cubes $\mu(Q_{x_i,k_i})$, (a)
follows.

Let us prove (b) now. We have
\begin{eqnarray*}
\sum_k \int |S_kf(x)|^2\,d\nu_k(x) & = & \int |S_k f(x)|^2\, d\nu(x,k) \\
& = & 2\int_0^\infty \lambda\,\nu\{(x,k):\, |S_kf(x)|>\lambda\}\,
d\lambda.
\end{eqnarray*}
We consider the maximal operator
$$M_S f(x) = \sup_{z,k:\,x\in (Q_{z,k})^\circ} |S_k f(z)|.$$
This operator is bounded on $L^2(\mu)$ (see Remark \ref{remms} below).
We set $E_\lambda = \{x\in\R^d:\, M_S f(x)>\lambda\}$. Then we have
$$\{(x,k):\, |S_kf(x)|>\lambda\}\subset \wt{E_\lambda},$$
By (a) we obtain
\begin{eqnarray*}
\sum_k \int |S_kf(x)|^2\,d\nu_k(x) & \leq &
C\int_0^\infty \lambda\,\nu(\wt{E_\lambda})\,d\lambda \\
& \leq & C\,C_9\int_0^\infty \lambda\,\mu(E_\lambda)\,d\lambda \\
& \leq & C\,C_9\int M_S f(x)^2\, d\mu(x) \,\leq\, C\,C_9\,\|f\|_{L^2(\mu)}^2,
\end{eqnarray*}
and we are done.
\end{proof}

\begin{rem}  \label{remms}
Consider the following non centered maximal operator
$$M_{(2)}f(x) = \sup_{Q:\,x\in Q} \frac{1}{\mu(2Q)} \int_Q |f(x)|\,d\mu(x).$$
This operator is bounded on $L^p(\mu)$, $1<p\leq\infty$, and of
weak type $(1,1)$ (see \cite[Section 6]{Tolsa3}).

It is not difficult to check that $M_Sf(x) \leq C\,M_{(2)}f(x)$ for all
$x\in\supp(\mu)$. Indeed, assume $x\in Q_{z,k}$ for some $z\in\supp(\mu)$,
$k\in\Z$, and let $N_0$ be the smallest integer such that
will be more involved than in \cite{DJS}, basically due to the fact t
$Q_{z,k-1}\subset 2^{N_0}Q_{z,k}$. We have
\begin{eqnarray*}
|S_kf(z)| & \leq & \left(\int_{Q_{z,k}}
+ \sum_{j=1}^{N_0}  \int_{2^jQ_{z,k}\setminus 2^{j-1}Q_{z,k}}\right)
s_k(z,y)\,|f(y)|\,d\mu(y) \\
& \leq &\frac{C}{\ell(Q_{z,k})^n}\int_{Q_{z,k}}|f|\,d\mu + \sum_{j=1}^{N_0}
\frac{C}{\ell(2^jQ_{z,k})^n} \int_{2^jQ_{z,k}}|f|\,d\mu \\
& \leq & C\,M_{(2)}f(x) + C\, \sum_{j=1}^{N_0}
\frac{\mu(2^{j+1}Q)}{\ell(2^{j+1}Q_{z,k})^n}\,M_{(2)}f(x) \\
& \leq & C\,(1+ \delta(Q_{z,k},\, 2^{N_0+1}Q_{z,k}))\,M_{(2)}f(x) \,\leq \,
C\,M_{(2)}f(x).
\end{eqnarray*}
\end{rem}

\vspace{3mm}
\begin{lemma} \label{parap}
If $g\in\rbmo(\mu)$ and $f\in L^2(\mu)$, then
$$\sum_{k\in\Z} \|(D_k^Ng)\cdot S_kf\|_{L^2(\mu)}^2 \leq C\,\|g\|_*^2\,
\|f\|_{L^2(\mu)}^2.$$
\end{lemma}

\begin{proof}
By Corollary \ref{corolprbmo} and the subsequent remark in \rf{lprbmo2}, since
$g\in \rbmo(\mu)$, we have
$$\sum_{j=k-2} \|D_j^N g\|_{L^2(\mu\mid Q)}^2 \leq C\,
\|g\|_*^2\,\mu(Q),$$
for any doubling cube $Q$ of the $k$th generation. Therefore, the lemma
follows from (b) in the preceding lemma taking $a_k:= (D_k^N\,g)^2$.
\end{proof}

As a direct consequence of the previous results we get:

\begin{lemma}
Given $b\in\rbmo(\mu)$,
the operators $U_{m,b}$ are bounded on $L^2(\mu)$ uniformly on $m$.
\end{lemma}

\begin{proof}
By Lemmas \ref{parap} and \ref{inverbmo}, for $f,g\in L^2(\mu)$ we have
\begin{eqnarray*}
|\langle U_{m,b}f,\, g\rangle| & \leq & \sum_k |\langle P_{k,b}^N S_k f,
\, D_k g\rangle| \\
& \leq  &
\left(\sum_k \|(D_k^N\Phi_N^{-1}(b))\cdot S_kf\|_{L^2(\mu)}^2\right)^{1/2}
\left( \sum_k \|D_kg\|_{L^2(\mu)}^2 \right)^{1/2} \\
& \leq & C \,\|f\|_{L^2(\mu)}\, \|\Phi_N^{-1}b\|_*\, \|g\|_{L^2(\mu)} \,
\leq \, C \,\|f\|_{L^2(\mu)}\, \|b\|_*\, \|g\|_{L^2(\mu)}.
\end{eqnarray*}
\end{proof}

In order to prove Theorem \ref{t1}, we want to apply the version of the $T(1)$
theorem in Lemma \ref{t10} to the operator
$T - U_{b_1}-  U_{b_2}^*$, with $b_1:= T(1)$ and $b_2:=T^*(1)$.
Given $b\in\rbmo(\mu)$,
we will not be able to show that $U_b$ is a CZO.
Instead we will show that $U_b$ satisfies some weaker assumptions, which will
be enough for our purposes.

\begin{lemma}  \label{weakczo}
Given $b\in\rbmo(\mu)$,
there are constants $C_{10},\, C_{11}$ such that,
for each $m$, the kernel $u^m(x,y)$ of $U_{m,b}$
satisfies the following properties:
\begin{itemize}
\item[(1)] $\ds |u^m(x,y)| \leq \frac{C_{10}}{|x-y|^n}$ if $x\neq y$.
\item[(2)] Let $x,x'\in\supp(\mu)$ with $x'\in Q_{x,h}$. Let $y\in\supp(\mu)$
be such that $y\in Q_{x,j}\setminus Q_{x,j+1}$ for some $j\leq h-10$.
Then,
$$|u^m(x,y) - u^m(x',y)| + |u^m(x,y) - u^m(x',y)|
\leq C_{11}\,\frac{|x-x'|}{\ell(Q_{x,j+4})\, |x-y|^n}.$$
\end{itemize}
\end{lemma}

Let us notice that the constants $C_{10},\,C_{11}$ above are independent of
$m$.

\begin{rem}
It is clear that any CZO satisfies the properties stated in the lemma above.
On the other hand, it can be seen that any operator fulfilling these properties
is a HCZO. Indeed, given $x,x'\in\supp(\mu)$ such that $x'\in Q_{x,h}\setminus
Q_{x,h+1}$, we have
\begin{eqnarray*}
\lefteqn{\int_{|x-y|\geq 2|x-x'|} |u^m(x,y) - u^m(x',y)|\,d\mu(y)} && \\
& \leq &
\int_{Q_{x,h-10}\setminus Q_{x,h+1}} \frac{C}{|x-y|^n}\, d\mu(y) \\
& & \mbox{} + C\,\sum_{i=10}^\infty \int_{Q_{x,h-i-1} \setminus Q_{x,h-i}}
\frac{|x-x'|}{\ell(Q_{x,h-i+4})\, |x-y|^n}\, d\mu(y) \\
& \leq & C +  C\,\sum_{i=10}^\infty \frac{|x-x'|}{\ell(Q_{x,h-i+4})}
\,\leq \,C.
\end{eqnarray*}
\end{rem}

\begin{proof}[Proof of Lemma \ref{weakczo}]
Let $u_k(x,y)$ be the kernel of $D_k\,P_{k,b}^N\,S_k$. Observe that
$$u_k(x,y) = \int d_k(x,z)\, a_k(z)\, s_k(z,y)\, d\mu(z),$$
where $a_k = D_k^N\,\Phi_N^{-1}(b)$. Since $\Phi_N^{-1}(b)\in\rbmo(\mu)$, it
follows that $a_k\in L^\infty(\mu)$, with $\|a_k\|_{L^\infty(\mu)} \leq
C\, \|b\|_*$ (the details are left to the reader).

Let us see that $u^m(x,y)$ satisfies condition (1). Take $x,y\in\supp(\mu)$ and
let $j\in\Z$ be such that $y\in Q_{x,j}\setminus Q_{x,j+1}$. We have
\begin{eqnarray*}
|u_k(x,y)|\! & \leq &\! C\int |d_k(x,z)\, s_k(z,y)|\, d\mu(z) \\
& = & \!C\int_{|x-z|\leq |x-y|/2} +\, C \int_{|x-z|> |x-y|/2}
\, \leq \,\,  \frac{C}{(\ell(Q_{x,k}) + |x-y|)^n}.
\end{eqnarray*}
Let us remark that the constant $C$ above equals $C\,\|b\|_*$.
It is not difficult to check that $\supp\, u_k(x,\cdot)\subset Q_{x,k-3}$.
Thus we have
$$|u^m(x,y)| \leq \sum_{k\leq j+3}  \frac{C}{(\ell(Q_{x,k}) + |x-y|)^n}
\leq \frac{C}{|x-y|^n}.$$

Now we will show that condition (2) of the lemma is also satisfied.
First we will deal with the term $|u^m(x,y) - u^m(x',y)|$. We have
$$I := |u_k(x,y) - u_k(x',y)| \leq C \int |(d_k(x,z)- d_k(x',z))\, s_k(z,y)|
\, d\mu(z).$$
We only have to estimate $I$ for $k-4\leq j$, because otherwise
$$\supp(u_k(x,\cdot) - u_k(x',\cdot)) \subset Q_{x,k-3} \cap Q_{x',k-3} \subset
Q_{x,j+1}.$$

Since $x'\in Q_{x,h}$ and $h\geq j+10>k$ (we are assuming $k-4\leq j$),
we have
$$|d_k(x,z) - d_k(x',z)| \leq C\,\frac{|x-x'|}{\ell(Q_{x,k})}\cdot
\frac{1}{(\ell(Q_{x,k}) + |x-z|)^n}.$$
Therefore,
\begin{eqnarray*}
I & \leq &  C\,\frac{|x-x'|}{\ell(Q_{x,k})} \int_{z\in Q_{x,k-3}}
\frac{1}{(\ell(Q_{x,k}) + |x-z|)^n}\,|s_k(z,y)| d\mu(z)\\
& = &  C\,\frac{|x-x'|}{\ell(Q_{x,k})}
\left( \int_{\begin{subarray}{l} z\in Q_{x,k-3}\\|x-z|\leq |x-y|/2
\end{subarray}} +
\int_{\begin{subarray}{l} z\in Q_{x,k-3}\\|x-z|> |x-y|/2\end{subarray}}\right)
\\ & \leq &  C\,\frac{|x-x'|}{\ell(Q_{x,k})}
\cdot \frac{1}{(\ell(Q_{x,k}) + |x-y|)^n}.
\end{eqnarray*}
We derive
\begin{eqnarray*}
|u^m(x,y) - u^m(x',y)| & \leq & C \sum_{k\leq j+4}
\frac{|x-x'|}{\ell(Q_{x,k})}
\cdot \frac{1}{(\ell(Q_{x,k}) + |x-y|)^n}\\
& \leq &
C\,\frac{|x-x'|}{\ell(Q_{x,j+4})\,|x-y|^n}.
\end{eqnarray*}

The estimates for the term $|u^m(y,x) - u^m(y,x')|$ are similar.
\end{proof}

\begin{rem}  \label{remlio}
The proof of Lemma \ref{t10}, corresponding to the $T(1)$ theorem in the case
$T(1)=T^*(1)=0$, given in the preceding
section can also be extended to operators
satisfying the assumptions of Lemma \ref{weakczo} (still assuming that the
kernel of the operator is bounded on $L^\infty$, for technical reasons).
This is an exercise which, again, is left to the reader.

Let us also notice that (as far as we know)
the arguments of Section \ref{sect10}
cannot be extended to HCZO's. Recall that even in the doubling situation
it is not known if the $T(1)$ theorem holds for HCZO's.
\end{rem}

Now we can finish the proof of the $T(1)$ theorem in the general case.

\begin{proof}[Proof of Theorem \ref{t1}]
We only have to prove that the operators $\wt{T}_\ve$ are bounded uniformly on
$\ve>0$. For a fixed $\ve>0$, we denote $b_1:= \wt{T}_\ve(1)$ and $b_2:=
\wt{T}^*_\ve(1)$. Since $b_1,b_2\in\bmo_\rho(\mu)$
and $\wt{T}_\ve$ is weakly bounded, from Lemma \ref{caracterist} it follows
$b_1,b_2\in \rbmo(\mu)$.
It is straightforward to check that $U_{b_1}^*(1) = U_{b_2}^*(1) =0$. We also
have $U_{b_i}(1)=b_i$, $i=1,2$. Indeed,
$$U_{m,b_i}(1) = \sum_{k=-m}^m D_k\, P_{k,b_i}^N\,S_k(1) =
 \sum_{k=-m}^m D_k\, D_k^N\, \Phi_N^{-1}(b_i).$$
Because of Lemma \ref{lemadk},
the operator $\sum_{k=-m}^m D_k D_k^N$ converges
strongly to $\Phi_N$ in $L^2(\mu)$ as $m\to+\infty$.
With estimates analogous to the ones
in Lemma \ref{weakczo}, it can be seen that $\sum_{k=-m}^m D_k D_k^N$ is
a HCZO, and taking into account
that $\sum_{k=-m}^m D_k D_k^N(1) = 0$, it follows that
this operator is bounded on $\rbmo(\mu)$ uniformly on $m$.
Arguing as in \cite[Lemma 2.9]{DJS}, it can be shown that for any function $g$
bounded with compact support and any $f\in\rbmo(\mu)$,
$$\lim_{m\to+\infty} \Bigl\langle \sum_{k=-m}^m D_k D_k^N (f),\,g\Bigr\rangle =
\bigl\langle \Phi_N(f),\,g\bigr\rangle.$$
As in \cite{DJS}, this implies $U_{b_i} (1) = b_i.$

Now it only remains to apply the version of the $T(1)$ theorem stated in Lemma
\ref{t10} to the operator $\wt{T}_\ve - U_{b_1} - U_{b_2}^*$.
Recall that that lemma applies to CZO's with bounded kernel. However, as
explained in the remark above, it is enough that the operator fulfil
condition (2) of Lemma \ref{weakczo}, instead of the usual gradient condition
demanded from the kernels of CZO's.
Moreover, the additional hypothesis in Lemma \ref{t10}
about the $L^\infty$-boundedness of the kernel was useful in the preceding
section to deal with the convergence of some integrals and also for
the proof Lemma \ref{convernormal}.
Although the kernels of $U_{b_1}$ and $U_{b_2}^*$ are not $L^\infty$-bounded,
we already know that these operators are bounded on $L^2(\mu)$ and that they
are weak limits of operators $U_{m,b_1}$, $U_{m,b_2}$, with `nice' kernels.
This allows the extension of the arguments for proving Lemma \ref{t10}
to the present situation. We omit the detailed arguments.
\end{proof}



\end{document}